\input amstex
\input amsppt.sty
\magnification=\magstep1
\hsize=30truecc
\vsize=22.2truecm
\baselineskip=16truept
\TagsOnRight
\nologo
\pageno=1
\topmatter
\def\N{\Bbb N}
\def\Z{\Bbb Z}
\def\Q{\Bbb Q}

\def\l{\left}
\def\r{\right}
\def\b{\bigg}

\def\({\b(}
\def\[{\b[}
\def\){\b)}
\def\]{\b]}

\def\t{\text}
\def\f{\frac}
\def\mo{\roman{mod}}
\def\ord{\roman{ord}}

\def\sm{\setminus}

\def\bi{\binom}
\def\eq{\equiv}

\def\ls{\leqslant}
\def\gs{\geqslant}
\def\al{\alpha}

\def\da{\delta}

\def\ms{\medskip}
\def\Proof{\noindent{\it Proof}}
\def\Remark{\medskip\noindent{\it Remark}}
\def\Ack{\noindent {\bf Acknowledgment}}

\def\st#1#2{\thickfrac\thickness0{#1}{#2}}
\def\nh{\noalign{\hrule}}
\def\hh{height4pt}

\def\om{&\omit &}

\hbox {Trans. Amer. Math. Soc. 359(2007), no.\,11, 5525--5553.}
\bigskip
\title Combinatorial congruences modulo prime powers\endtitle
\author Zhi-Wei Sun$^1$ and Donald M. Davis$^2$\endauthor
\leftheadtext{Zhi-Wei Sun and Donald M. Davis}
\affil $^1$Department of Mathematics, Nanjing University
\\ Nanjing 210093, People's Republic of China
\\zwsun\@nju.edu.cn
\\ {\tt http://math.nju.edu.cn/$^\sim$zwsun}
\medskip
$^2$Department of Mathematics, Lehigh University\\Bethlehem, PA
18015, USA
\\dmd1\@lehigh.edu
\\ {\tt http://www.lehigh.edu/$^{\sim}$dmd1}
\endaffil
\abstract
Let $p$ be any prime,
and let $\al$ and $n$ be nonnegative integers.
Let $r\in\Z$ and $f(x)\in\Z[x]$.
We establish the congruence
$$p^{\deg f}\sum_{k\eq r\, (\mo\ p^{\al})}\bi nk(-1)^kf\l(\f{k-r}{p^{\al}}\r)
\eq0\ \l(\mo\ p^{\sum_{i=\al}^{\infty}\lfloor  n/{p^i}\rfloor}\r)$$
(motivated by a conjecture
arising from algebraic topology),
and obtain the following vast generalization of Lucas' theorem:
If $\al$ is greater than one, and $l,s,t$ are nonnegative integers with $s,t<p$,
then
$$\aligned&\f1{\lfloor n/p^{\al-1}\rfloor!}
\sum_{k\eq r\,(\mo\ p^{\al})}\bi{pn+s}{pk+t}(-1)^{pk}\l(\f{k-r}{p^{\al-1}}\r)^l
\\\eq& \f1{\lfloor n/p^{\al-1}\rfloor!}
\sum_{k\eq r\,(\mo\ p^{\al})}\bi nk\bi st(-1)^k\l(\f{k-r}{p^{\al-1}}\r)^l
\ (\mo\ p).
\endaligned$$
We also present an application of the first congruence to Bernoulli polynomials,
and apply the second congruence to show that a $p$-adic order bound given
by the authors in a previous
paper can be attained when $p=2$.
\endabstract
\thanks 2000 {\it Mathematics Subject Classification}.\,Primary 11B65;
Secondary 05A10, 11A07, 11B68, 11S05.\newline\indent The first
author is responsible for communications, and supported
by the National Science Fund for Distinguished Young Scholars in P. R. China
(Grant No. 10425103).
\endthanks
\endtopmatter
\document

\heading{1. Introduction}\endheading

In this paper we establish a number of new congruences for sums
involving binomial coefficients with the summation index
restricted in a residue class modulo a prime power; some of them
are vast extensions of some classical congruences. We begin by
providing some historical background for these results.

Let $p$ be a prime, and let $\Q_p$ and $\Z_p$ denote the field of $p$-adic numbers
and the ring of $p$-adic integers respectively.
For $\omega\in\Q_p\sm\{0\}$ we define its {\it $p$-adic order} by
$\ord_p(\omega)=\max\{a\in\Z:\,\omega/p^a\in\Z_p\}$; in addition, we set
$\ord_p(0)=+\infty$.

In 1913, A. Fleck (cf. [D, p.\,274]) proved that for any
$n\in\Z^+=\{1,2,3,\ldots\}$ and $r\in\Z$
we have the congruence
$$\sum_{k\eq r\,(\mo\ p)}\bi nk(-1)^k\eq0\ \ \l(\mo\ p^{\lfloor\f{n-1}{p-1}\rfloor}\r)$$
(where $\lfloor\cdot\rfloor$
is the greatest integer function, and we regard $\bi xk=0$ for $k=-1,-2,-3,\ldots$);
that is,
$$\ord_p\(\sum_{k\eq r\,(\mo\ p)}\bi nk(-1)^k\)\gs\l\lfloor\f{n-1}{p-1}\r\rfloor.$$
In 1977, C. S. Weisman [We] showed further that
if $\al,n\in\N=\{0,1,2,\ldots\}$ and $r\in\Z$ then
$$\ord_p\(\sum_{k\eq r\,(\mo\ p^{\al})}\bi nk(-1)^k\)
\gs\l\lfloor\f{n-p^{\al-1}}{\varphi(p^{\al})}\r\rfloor,$$
where $\varphi$ is the well-known Euler function.
Weisman remarked that this kind of work is closely related to $p$-adic continuation.

In 2005, motivated by Fontaine's theory of
$(\phi,\Gamma)$-modules, D. Wan got an extension of Fleck's result
in his lecture notes by giving a lower bound for the $p$-adic
order of the sum $\sum_{k\eq r\,(\mo\ p)}\bi
nk(-1)^k\bi{(k-r)/p}l,$ where $l,n\in\N$ and $r\in\Z$. Soon after
this, Z. W. Sun [S06] obtained a common generalization of
Weisman's and Wan's extensions of Fleck's congruence by studying
the $p$-adic order of the sum
$$\sum_{k\eq r\,(\mo\ p^{\al})}\bi nk(-1)^k\bi{(k-r)/p^{\al}}l$$
(with $\al\in\N$) via a combinatorial approach.
This direction is closely related to the $\al$th iteration of the basic $\psi$-operator
in Iwasawa theory and $p$-adic Langlands correspondence (cf. P. Colmez [C]
and Wan [W]).
However, when $l\gs n/p^{\al}$,
any result along this line
yields no nontrivial lower bound for the $p$-adic order of the last sum.

Unlike the above development of Fleck's congruence,
there is another direction motivated by algebraic topology.
In order to obtain a strong lower bound for
homotopy exponents
of the special unitary group
$\t {SU}(n)$, the authors [DS]
were led to show that if $\al,n\in\N$ and $r\in\Z$ then
$$\min_{f(x)\in\Z[x]}\ord_p\(\sum_{k\eq r\,(\mo\ p^{\al})}
\bi nk(-1)^kf\l(\f{k-r}{p^{\al}}\r)\)
\gs\ord_p\l(\l\lfloor\f{n}{p^{\al}}\r\rfloor!\r).$$
As shown in [DS], this inequality
implies a subtle divisibility property of Stirling numbers of the second kind.
Note that if $n\in\Z^+$ then
$$\ord_p(n!)=\sum_{i=1}^{\infty}\l\lfloor \f n{p^i}\r\rfloor
<\sum_{i=1}^{\infty}\f n{p^i}=\f n{p(1-p^{-1})}=\f n{p-1}$$
and hence $\ord_p(n!)\ls (n-1)/(p-1)$.

Now we introduce some conventions used throughout this paper.
As usual, the degree of the zero polynomial
is regarded as $-\infty$.
For an integer $a$ and a positive real number $m$, we let
$\{a\}_m$ denote the fractional part of $a/m$ times $m$
(i.e., $\{a\}_m$ is the unique number in the interval $[0,m)$
with $a-\{a\}_m\in m\Z$). For a prime $p$, if $a,b\in\N$, then
$\tau_p(a,b)$ stands for the number of carries when adding
$a$ and $b$ in base $p$; a theorem of E. Kummer states that
$\tau_p(a,b)=\ord_p\bi{a+b}a$.
\smallskip

Here is our first theorem.

\proclaim{Theorem 1.1} Let $p$ be a prime and $f(x)\in\Z_p[x]$.
Let $\al,n\in\N$ and $r\in\Z$.
Then
$$\aligned&p^{\deg f}\sum_{k\eq r\,(\mo\ p^{\al})}\bi nk(-1)^kf\l(\f{k-r}{p^{\al}}\r)
\\\eq&0\ \l(\mo\ p^{\sum_{i=\al}^{\infty}\lfloor n/p^{i}\rfloor
+\tau_p(\{r\}_{p^{\al-1}},\{n-r\}_{p^{\al-1}})}\r),
\endaligned$$
i.e.,
$$\aligned&\ord_p\(\sum_{k\eq r\,(\mo\ p^{\al})}\bi nk(-1)^kf\l(\f{k-r}{p^{\al}}\r)\)
\\\gs&\ord_p\l(\l\lfloor\f n{p^{\al-1}}\r\rfloor!\r)-\deg f
+\tau_p(\{r\}_{p^{\al-1}},\{n-r\}_{p^{\al-1}}).
\endaligned\tag1.1 $$
\endproclaim
\Remark\ 1.1. It is interesting to compare (1.1)
with the following inequality ([Theorem 5.1, DS])
established for a certain topological purpose:
$$\aligned&\ord_p\(\sum_{k\eq r\,(\mo\ p^{\al})}\bi nk(-1)^kf\l(\f{k-r}{p^{\al}}\r)\)
\\&\ \gs\ord_p\l(\l\lfloor\f n{p^{\al}}\r\rfloor!\r)
+\tau_p(\{r\}_{p^{\al}},\{n-r\}_{p^{\al}}).\endaligned\tag1.2 $$
Note that $\ord_p(\lfloor n/p^{\al}\rfloor!)
=\ord_p(\lfloor n/p^{\al-1}\rfloor!)-\lfloor n/p^{\al}\rfloor$ and
$$0\ls\tau_p(\{r\}_{p^{\al}},\{n-r\}_{p^{\al}})
-\tau_p(\{r\}_{p^{\al-1}},\{n-r\}_{p^{\al-1}})\ls 1.$$
In general, when $\deg f<\lfloor n/p^{\al}\rfloor$, the term
$\tau_p(\{r\}_{p^{\al-1}},\{n-r\}_{p^{\al-1}})$ in (1.1) cannot be replaced by
$\tau_p(\{r\}_{p^{\al}},\{n-r\}_{p^{\al}})$.
A principal motivation for the development of Theorem 1.1 is
that the inequality $(1.2)$,
although quite sharp when $\deg f\gs\lfloor n/p^\al\rfloor$,
is not a good estimate for smaller values of $\deg f$.
\medskip

{\it Example} 1.1.  Let $p=\al=2$, $r=1$ and $n=20$.
Then, for each $0<l<\lfloor n/p^{\al}\rfloor=5$,
equality in (1.1) with $f(x)=x^l$ is attained while
$$\tau_p(\{r\}_{p^{\al}},\{n-r\}_{p^{\al}})
=\tau_2(1,3)=\tau_2(1,1)+1=\tau_p(\{r\}_{p^{\al-1}},\{n-r\}_{p^{\al-1}})+1.$$

{\it Example} 1.2. Let $p=3$, $\al=r=2$, $90\ls n\ls 98$
and $0\ls l<\lfloor n/p^{\al}\rfloor=10$.
In Table 1 below, $\da_n(l)$ denotes the left hand side of (1.1) minus the right hand side
with $f(x)=x^l$.

\bigskip
\centerline{Table 1: Values of $\da_n(l)$ with $0\ls l\ls 9$ and $90\ls n\ls 98$}
\smallskip
\centerline{\vbox{\offinterlineskip
\halign{\vrule#&\ \ #\ \hfill   &&\vrule#&\ \ \hfill#\ \
\cr\nh\cr \hh\om\om\om\om\om\om\om\om\om\om\om
\cr &$\st {\ \ l}{n\ \ }$& &0\,& &$1\,$& &$2\,$& &$3\,$& &$4\,$& &$5\,$& &$6\,$& &$7\,$& &$8\,$& &$9\,$&
\cr \hh\om\om\om\om\om\om\om\om\om\om\om
\cr\nh\cr \hh\om\om\om\om\om\om\om\om\om\om\om
\cr & $90$ && $\,1$ &&$\,0$ && $\,2$ && $\,0$ && $\,1$ && $\,0$ && $\,1$ &&$\,0$ && $\,3$ && $\,0$&
\cr \hh\om\om\om\om\om\om\om\om\om\om\om
\cr\nh\cr \hh\om\om\om\om\om\om\om\om\om\om\om
\cr & $91$ && $1$ && $0$ && $1$ && $0$ && $3$ && $0$ && $1$ &&$0$ && $1$ && $0$&
\cr \hh\om\om\om\om\om\om\om\om\om\om\om
\cr \nh\cr \hh\om\om\om\om\om\om\om\om\om\om\om
\cr & $92$ && $0$ && $1$ && $0$ && $3$ && $0$ && $1$ && $0$ &&$1$ && $0$ && $2$&
\cr \hh\om\om\om\om\om\om\om\om\om\om\om
\cr \nh\cr \hh\om\om\om\om\om\om\om\om\om\om\om
\cr & $93$ && $0$ && $2$ && $0$ && $1$ && $0$ && $1$ && $0$ &&$4$ && $0$ && $1$&
\cr \hh\om\om\om\om\om\om\om\om\om\om\om
\cr \nh\cr \hh\om\om\om\om\om\om\om\om\om\om\om
\cr & $94$ && $0$ && $1$ && $0$ && $2$ && $0$ && $1$ && $0$ &&$1$ && $0$ && $3$&
\cr \hh\om\om\om\om\om\om\om\om\om\om\om
\cr \nh\cr \hh\om\om\om\om\om\om\om\om\om\om\om
\cr & $95$ && $1$ && $0$ && $0$ && $0$ && $0$ && $1$ && $1$ &&$0$ && $0$ && $0$&
\cr \hh\om\om\om\om\om\om\om\om\om\om\om
\cr \nh\cr \hh\om\om\om\om\om\om\om\om\om\om\om
\cr & $96$ && $1$ && $0$ && $0$ && $0$ && $0$ && $1$ && $2$ &&$0$ && $0$ && $0$&
\cr \hh\om\om\om\om\om\om\om\om\om\om\om
\cr \nh\cr\hh\om\om\om\om\om\om\om\om\om\om\om
\cr & $97$ && $1$ && $0$ && $0$ && $0$ && $0$ && $1$ && $1$ &&$0$ && $0$ && $0$&
\cr \hh\om\om\om\om\om\om\om\om\om\om\om
\cr \nh\cr \hh\om\om\om\om\om\om\om\om\om\om\om
\cr & $98$ && $1$ && $3$ && $0$ && $1$ && $0$ && $1$ && $1$ &&$4$ && $0$ && $1$&
\cr \hh\om\om\om\om\om\om\om\om\om\om\om
\cr \nh\cr }}}

\bigskip
\bigskip
{\it Example} 1.3.  The combination of using Theorem
1.1 when $\deg f\ls \lfloor n/p^\al\rfloor$ and the inequality $(1.2)$
for larger values of $\deg f$ provides an excellent estimate for
$$\ord_p\(\sum_{k\eq r\,(\mo\ p^{\al})}
\bi nk(-1)^kf\l(\f{k-r}{p^{\al}}\r)\).$$
For example, if $p=2$, $\al=1$, $r=0$, $n=20$ and $f(x)=x^l$,
then the actual values of the expression for $l=0,\ldots,21$ are
$$19,19,17,17,14,14,12,12,10,10,8,8,8,8,11,9,9,9,8,8,8,8.$$
The inequality $(1.2)$ guarantees that each of the 22 numbers should be at least $8$,
while Theorem 1.1
gives the lower bounds
$$18, 17, 16, 15, 14, 13, 12, 11, 10, 9, 8$$
for the first 11 of the above 22 values respectively.
The bound given by (1.1) is attained in this example
with $l=4$, while
$\ord_p(\lfloor(n/p^{\al-1})\rfloor!)=18<19=\lfloor(n-p^{\al-1})/\varphi(p^{\al})\rfloor$;
this shows that for a general $f(x)\in\Z_p[x]$ we cannot replace
$\ord_p(\lfloor(n/p^{\al-1})\rfloor!)$ in (1.1) by
Weisman's bound $\lfloor(n-p^{\al-1})/\varphi(p^{\al})\rfloor$ even if $r=0$
(and hence the $\tau$-term vanishes).
\medskip

Here is an application of Theorem 1.1 to Bernoulli polynomials.
(The reader is referred to [IR] and [S03] for some basic properties and known congruences
concerning Bernoulli polynomials.)

\proclaim{Corollary 1.1} Let $p$ be a prime, and let $\al\in\N$, $m,n\in\Z^+$ and $r\in\Z$. Then
$$\aligned&\ord_p\(\f{p^{m-1}}m
\sum_{k=0}^n\bi nk(-1)^kB_m\(\l\lfloor\f{k-r}{p^{\al}}\r\rfloor\)\)
\\&\gs\sum_{i=\al}^{\infty}\l\lfloor\f{n-1}{p^i}\r\rfloor
+\tau_p(\{r-1\}_{p^{\al-1}},\{n-r\}_{p^{\al-1}}).
\endaligned\tag1.3$$
\endproclaim
\Proof. Set $\bar r=r+p^{\al}-1$. In view of [S06, Lemma 2.1] and
the known identity $B_m(x+1)-B_m(x)=mx^{m-1}$, we have
$$\align&\sum_{k=0}^n\bi nk\f{(-1)^k}mB_m\l(\l\lfloor\f{k-r}{p^{\al}}\r\rfloor\r)
\\=&\sum_{k\eq \bar r\,(\mo\ p^{\al})}\bi {n-1}k\f{(-1)^{k-1}}m
\(B_m\l(\f{k-\bar r}{p^{\al}}+1\r)-B_m\l(\f{k-\bar r}{p^{\al}}\r)\)
\\=&\sum_{k\eq \bar r\,(\mo\ p^{\al})}\bi {n-1}k(-1)^{k-1}
\l(\f{k-\bar r}{p^{\al}}\r)^{m-1}.
\endalign$$
This, together with Theorem 1.1, yields that
$$\align&\ord_p\(\sum_{k=0}^n\bi nk\f{(-1)^k}mB_m\l(\l\lfloor\f{k-r}{p^{\al}}\r\rfloor\r)\)
\\\gs&\ord_p\l(\l\lfloor\f{n-1}{p^{\al-1}}\r\rfloor!\r)-(m-1)
+\tau_p(\{\bar r\}_{p^{\al-1}},\{n-1-\bar r\}_{p^{\al-1}}).
\endalign$$
So (1.3) follows.  \qed

\bigskip
Let $p$ be any prime, $f(x)\in\Q_p[x]$ and $\deg f\ls l\in\N$.
It is well known that $f(a)\in\Z_p$ for all $a\in\Z$ if and only if
$f(x)=\sum_{j=0}^la_j\bi xj$ for some $a_0,\ldots,a_l\in\Z_p$.
Also, $f(x)\in\Z_p[x]$ if and only if it can be written
in the form $\sum_{j=0}^lb_j(x)_j$ with $b_j\in\Z_p$ and $(x)_j=j!\bi xj\in\Z[x]$.
(Recall that $x^l=\sum_{j=0}^lS(l,j)(x)_j$
where $S(l,j)\ (0\ls j\ls l)$ are Stirling numbers of the second kind.)
Thus, we can reformulate Theorem 1.1 as follows.

\proclaim{Theorem 1.2} Let $p$ be a prime, $\al,n\in\N$ and $r\in\Z$.
Let $f(x)\in\Q_p[x]$, $\deg f\ls l\in\N$, and $f(a)\in\Z_p$ for all $a\in\Z$. Then we have
$$\aligned&\ord_p\(\sum_{k\eq r\,(\mo\ p^{\al})}\bi nk(-1)^kf\l(\f{k-r}{p^{\al}}\r)\)
\\\gs&\ord_p\l(\l\lfloor\f n{p^{\al-1}}\r\rfloor!\r)-l-\ord_p(l!)
+\tau_p(\{r\}_{p^{\al-1}},\{n-r\}_{p^{\al-1}}).
\endaligned\tag1.4$$
\endproclaim

\Remark\ 1.2. This theorem has topological background.
In the case $p=\al=r=2$ and $f(x)=\bi xl$, it first arose as a conjecture of
the second author in his study of algebraic topology.
\medskip

Let $[x^n]F(x)$ denote the coefficient of $x^n$ in the power series expansion of $F(x)$.
Theorem 1.1 also has the following equivalent form.

\proclaim{Theorem 1.3} Let $p$ be a prime, and let
$\al,l,n,r\in\N$. If $r>n-(l+1)p^{\al}$, then
$$\aligned&\ord_p\([x^r]\f{(1-x)^n}{(1-x^{p^{\al}})^{l+1}}\)
\\\gs&\ord_p\l(\l\lfloor\f{n}{p^{\al-1}}\r\rfloor!\r)-l-\ord_p(l!)
+\tau_p(\{r\}_{p^{\al-1}},\{n-r\}_{p^{\al-1}}).
\endaligned\tag1.5$$
\endproclaim

\Proof. Let $r>n-(l+1)m$ where $m=p^{\al}$.
Observe that
$$\align [x^r]\f{(1-x)^n}{(1-x^m)^{l+1}}
=&\sum_{k=0}^r\bi nk(-1)^k[x^{r-k}](1-x^m)^{-l-1}
\\=&\sum_{k=0}^r\bi nk(-1)^k[x^{r-k}]\sum_{j\gs0}\bi{l+j}l(x^m)^j
\\=&\sum\Sb 0\ls k\ls r\\k\eq r\,(\mo\ m)\endSb\bi nk(-1)^k\bi{l-(k-r)/m}l.
\endalign$$
If $r<k\ls n$ and $k\eq r\ (\mo\ m)$, then
$0<(k-r)/m\ls (n-r)/m<l+1$ and hence $\bi{l-(k-r)/m}l=0$.
Therefore
$$\align[x^r]\f{(1-x)^n}{(1-x^m)^{l+1}}
=&\sum\Sb 0\ls k\ls n\\k\eq r(\mo\ m)\endSb\bi nk(-1)^k\bi{l-(k-r)/m}l,
\\=&(-1)^l\sum_{k\eq r'\,(\mo\ m)}\bi nk(-1)^k\bi{(k-r')/m}l,
\endalign$$
where $r'=r+m$. Applying Theorem 1.1 or 1.2 we immediately get the inequality (1.5).
\qed

\smallskip

Here is another equivalent version of Theorem 1.1.

\proclaim{Theorem 1.4} Let $p$ be a prime and $\al$ be a nonnegative integer.
Let $r\in\Z$, and $f(x)\in\Z_p[x]$ with $\deg f=l\in\N$. Then,
there is a sequence $\{a_k\}_{k\in\N}$ of $p$-adic integers such
that for any $n\in\N$ we have
$$\aligned&\sum_{k=0}^n\bi nk(-1)^k\l\lfloor\f k{p^{\al-1}}\r\rfloor!
\bi{\{r\}_{p^{\al-1}}+\{k-r\}_{p^{\al-1}}}
{\{r\}_{p^{\al-1}}}a_k
\\&\ \ =\cases p^lf(\f{n-r}{p^{\al}})&\t{if}\ n\eq r\,(\mo\ p^{\al}),
\\0&\t{otherwise}.\endcases
\endaligned\tag1.6$$
\endproclaim
\Proof. The binomial inversion formula (cf. [GKP, (5.48)]) states that
$\sum_{k=0}^n\bi nk(-1)^kb_k=d_n$ for all $n\in\N$, if and only if
$\sum_{k=0}^n\bi nk(-1)^kd_k=b_n$ for all $n\in\N$. Thus the desired result
has the following equivalent form: There exists a sequence $\{a_n\}_{n\in\N}$
of $p$-adic integers such that for all $n\in\N$ we have
$$\align&\l\lfloor\f n{p^{\al-1}}\r\rfloor!\bi{\{r\}_{p^{\al-1}}+\{n-r\}_{p^{\al-1}}}
{\{r\}_{p^{\al-1}}}a_n
\\&\ =\sum_{k\eq r\,(\mo\ p^{\al})}\bi nk(-1)^k
p^lf\l(\f{k-r}{p^{\al}}\r).\endalign$$
This is essentially what Theorem 1.1 says. \qed

A famous theorem of E. Lucas states that
if $p$ is a prime and $n,r,s,t$ are nonnegative integers with $s,t<p$ then
$$\bi{pn+s}{pr+t}\eq\bi nr\bi st\ (\mo\ p).$$
Now we present our following analogue of Lucas' theorem.

\proclaim{Theorem 1.5} Let $p$ be a prime and $\al\gs2$ be an integer.
Then, for any $l,n\in\N$ and $r\in\Z$, we have the congruence
$$T_{l,\al+1}^{(p)}(n,r)\eq(-1)^{\{r\}_p}\bi{\{n\}_p}{\{r\}_p}
T_{l,\al}^{(p)}\l(\l\lfloor \f np\r\rfloor,
\l\lfloor\f rp\r\rfloor\r)\ (\mo\ p),\tag1.7$$
where
$$T_{l,\al}^{(p)}(n,r):=\f{l!p^l}{\lfloor n/p^{\al-1}\rfloor!}
\sum_{k\eq r\,(\mo\ p^{\al})}\bi nk(-1)^k\bi{(k-r)/p^{\al}}l.$$
\endproclaim

\Remark\ 1.3. Theorem 1.2 guarantees that $T_{l,\al}^{(p)}(n,r)\in\Z_p$
(and our proof of Theorem 1.2 given later is based on analysis of this $T$).
Theorem 1.5 provides further information on $T_{l,\al}^{(p)}(n,r)$ modulo $p$.
\medskip

Since $f(x)=l!\bi xl-x^l\in\Z[x]$ has degree smaller than $l$, we have
$$T_{l,\al}^{(p)}(n,r)\eq\f{p^l}{\lfloor n/p^{\al-1}\rfloor!}
\sum_{k\eq r\,(\mo\ p^{\al})}\bi nk(-1)^k\l(\f{k-r}{p^{\al}}\r)^l\ (\mo\ p)$$
by Theorem 1.1. Thus, Theorem 1.5 has the following equivalent version.

\proclaim{Theorem 1.6} Let $p$ be any prime,
and let $l,n\in\N$, $r,s,t\in\Z$ and $0\ls s,t<p$.
Then, for every $\al=2,3,\ldots$, we have
$$\aligned&\f1{\lfloor n/p^{\al-1}\rfloor!}
\sum_{k\eq r\,(\mo\ p^{\al})}\bi{pn+s}{pk+t}(-1)^{pk}\l(\f{k-r}{p^{\al-1}}\r)^l
\\\eq& \f1{\lfloor n/p^{\al-1}\rfloor!}
\sum_{k\eq r\,(\mo\ p^{\al})}\bi nk\bi st(-1)^k\l(\f{k-r}{p^{\al-1}}\r)^l
\ (\mo\ p).
\endaligned\tag1.8 $$
\endproclaim

\Remark\ 1.4. Theorem 1.6 is a vast generalization of Lucas' theorem.
Given a prime $p$ and nonnegative integers $n,r,s,t$ with $r\ls n$ and $s,t<p$,
if we apply (1.8) with $\al>\log_p(\max\{n,p\})$ and $l=0$
then we obtain Lucas' congruence
$\bi{pn+s}{pr+t}\eq\bi nr\bi st\ (\mo\ p)$.
We conjecture that (1.8) (or its equivalent form (1.7)) also holds when $\al=1$.
\medskip

Let $p>3$ be a prime. A well-known theorem of Wolstenholme asserts that
$\bi{2p-1}{p-1}\eq 1\ (\mo\ p^3)$, i.e., $\bi{2p}p\eq 2\ (\mo\ p^3)$.
In 1952 W. Ljunggren (cf. [B] [G]) generalized this as follows:
$\bi{pn}{pr}\eq\bi nr\ (\mo\ p^3)$
for any $n,r\in\N$. Our following conjecture extends Ljunggren's result greatly.

\proclaim{Conjecture 1.1} Let $p$ be an odd prime,
and let $\al\in\Z^+$, $n\in\N$ and $r\in\Z$.
Then, for all $l\in\N$ we have
$$T_{l,\al+1}^{(p)}(pn,pr)-T_{l,\al}^{(p)}(n,r)\eq\ \cases0\ (\mo\ p^3)&\t{if}\ p>3,
\\0\ (\mo\ p^2)&\t{if}\ p=3.\endcases$$
Equivalently, for any $f(x)\in\Z_p[x]$ we have
$$\align&\f1{\lfloor n/p^{\al-1}\rfloor!}
\sum_{k\eq r\,(\mo\ p^{\al})}\(\bi{pn}{pk}-\bi nk\)(-1)^k
f\l(\f{k-r}{p^{\al-1}}\r)
\\&\qquad\quad\eq\ \cases0\ (\mo\ p^3)&\t{if}\ p>3,
\\0\ (\mo\ p^2)&\t{if}\ p=3.\endcases
\endalign$$
\endproclaim

\Remark\ 1.5. The reason for the equivalence of the two parts in Conjecture 1.1
is as follows: For any $l\in\N$ we have
$$p^lx^l=\sum_{j=0}^l S(l,j)p^{l-j}j!p^j\bi xj
\ \t{and}\ l!p^l\bi xl=\sum_{j=0}^l(-1)^{l-j}s(l,j)p^{l-j}(p^jx^j),$$
where $s(l,j)\ (0\ls j\ls l)$ are Stirling numbers of the first kind.
\medskip

The proof of Theorem 1.5 involves our following refinement of Weisman's result.

\proclaim{Theorem 1.7} Let $p$ be any prime, and let $\al,n\in\N$,
$\al\gs2$, $r,s,t\in\Z$ and $0\ls s,t<p^{\al-2}$. Then
$$\aligned &p^{-\l\lfloor\f{p^{\al-2}n+s-p^{\al-1}}
{\varphi(p^{\al})}\r\rfloor}\sum_{k\eq p^{\al-2}r+t\,(\mo\ p^{\al})}
\bi {p^{\al-2}n+s}k(-1)^k
\\&\eq(-1)^t\bi{s}t\(p^{-\l\lfloor \f{n-p}{\varphi(p^2)}\r\rfloor}
\sum_{k\eq r\,(\mo\ p^2)}\bi{n}k(-1)^k\)\ (\mo\ p).
\endaligned\tag1.9$$
\endproclaim

Here is a consequence of this theorem.
\proclaim{Corollary 1.2} Let $\al\in\{2,3,\ldots\}$, $n\in\N$ and $r\in\Z$. Then
$$2^{-\l\lfloor\f{n-2^{\al-1}}{\varphi(2^{\al})}\r\rfloor}\sum_{k\eq r\,(\mo\ 2^{\al})}
\bi nk\eq1\ (\mo\ 2)$$
if and only if
$ \bi{\{n\}_{2^{\al-2}}}{\{r\}_{2^{\al-2}}}\eq1\ (\mo\ 2)$
and the following $(a)$ or $(b)$ holds.

{\rm (a)} $n_*>2$ and $n_*\not\eq2r_*+2\ (\mo\ 4)$,
\quad {\rm (b)} $n_*=2$ and $2\mid r_*$,
\newline
where $n_*=\lfloor n/2^{\al-2}\rfloor$ and $r_*=\lfloor r/2^{\al-2}\rfloor$.
\endproclaim

As a complement to Theorem 1.7, we have the following conjecture.

\proclaim{Conjecture 1.2} Let $p$ be any prime,
and let $n\in\N$, $r\in\Z$ and $s\in\{0,\ldots,p-1\}$.
If $p\mid n$ or $p-1\nmid n-1$, then
$$\aligned&p^{-\l\lfloor\f{pn+s-p}{\varphi(p^2)}\r\rfloor}
\sum_{k\eq pr+t\,(\mo\ p^2)}\bi{pn+s}k(-1)^k
\\\eq&(-1)^t\bi stp^{-\lfloor\f{n-1}{p-1}\rfloor}\sum_{k\eq r\,(\mo\ p)}\bi nk(-1)^k
\ \ (\mo\ p)\endaligned$$
for every $t=0,\ldots,p-1$.
When $s\not=p-1$, $p\nmid n$ and $p-1\mid n-1$, the least nonnegative residue of
$$p^{-\l\lfloor\f{pn+s-p}{\varphi(p^2)}\r\rfloor}
\sum_{k\eq pr+t\,(\mo\ p^2)}\bi{pn+s}k(-1)^k$$
modulo $p$ does not depend on $r$ for $t=s+1,\ldots,p-1$; moreover these residues
form a permutation of $1,\ldots,p-1$ if $s=0$ and $n\not=1$.
\endproclaim

Conjecture 5.2 of [DS] states that the bound in the inequality $(1.2)$
is attained if $f(x)=x^l$
and $l$ satisfies a certain congruence equation. In the following result, we prove this
conjecture in the case $p=2$ and $r=0$.

\proclaim{Theorem 1.8} Let $\al,n\in\N$, $l\gs\lfloor n/2^\al\rfloor\gs1$
and $$l\equiv\l\lfloor \f n{2^\al}\r\rfloor\ \l(\mo\ 2^{\lfloor\log_2(n/2^\al)\rfloor}\r).$$
Then
$$\ord_2\(\sum_{k\eq 0\,(\mo\ 2^{\al})}\binom nk(-1)^k\(\frac k{2^\al}\)^l\)
=\ord_2\l(\l\lfloor\f{n}{2^{\al}}\r\rfloor!\r).\tag1.10$$
\endproclaim

\Remark\ 1.6. Theorem 1.8 in the case $\al=0$ essentially asserts that
if $n\in\Z^+$ and $q\in\N$ then
$S(n+2^{\lfloor \log_2n\rfloor}q,\,n)$ is odd. This is because
$\sum_{k=0}^n\bi nk(-1)^kk^l=(-1)^nn!S(l,n)$ for $l\in\N$ (cf. [LW, pp.\,125--126]).
\medskip

The following conjecture is a refinement of [DS, Conjecture 5.2].

\proclaim{Conjecture 1.3} Let $p$ be a prime and $r$ be an integer.
Let $\al,n\in\N$, $n\gs2p^{\al}-1$, $l\gs\lfloor n/p^{\al}\rfloor$ and
$$l\eq\l\lfloor\f r{p^{\al}}\r\rfloor+\l\lfloor\f{n-r}{p^{\al}}\r\rfloor
\ \l(\mo\ (p-1)p^{\lfloor\log_p(n/p^{\al})\rfloor}\r).$$
Then
$$\f1{\lfloor n/p^{\al}\rfloor!\bi{n_*}{r_*}}
\sum_{k\eq r\,(\mo\ p^{\al})}\bi nk(-1)^k\l(\f{k-r}{p^{\al}}\r)^{l}
\eq\pm 1\ \ (\mo\ p),$$
where $r_*=\{r\}_{p^{\al}}$ and $n_*=r_*+\{n-r\}_{p^{\al}}$.
\endproclaim

For convenience, throughout this paper
we use $[\![A]\!]$ to denote the characteristic function of an assertion $A$, i.e.,
the value of $[\![A]\!]$ is $1$ or $0$ according to whether $A$ holds or not.
For $m,n\in\N$ the Kronecker symbol $\da_{m,n}$ stands for $[\![m=n]\!]$.

The next section is devoted to the proof of an equivalent version of Theorems 1.1--1.4.
In Section 3 we will prove Theorem 1.7 and Corollary 1.2
on the basis of Theorem 3.1,
and in Section 4 we establish Theorems 1.5 and 1.8.

\heading{2. Proof of Theorem 1.2}\endheading

In this section we prove the following equivalent version of Theorem 1.2.

\proclaim{Theorem 2.1} Let $p$ be a prime, and let $\al,l,n\in\N$.
Then, for all $r\in\Z$, we have
$$\ord_p(T_l(n,r))\gs\tau_p(\{r\}_{p^{\al-1}},\{n-r\}_{p^{\al-1}}),\tag2.1$$
where $T_l(n,r)$ stands for $T_{l,\al}^{(p)}(n,r)$ given in Theorem 1.5.
\endproclaim

\proclaim{Lemma 2.1} Theorem 2.1 holds in the case $\al=0$.
\endproclaim
\Proof. Clearly $\tau_p(\{r\}_{p^{-1}},\{n-r\}_{p^{-1}})=\tau_p(0,0)=0$.
Provided $\al=0$ we have
$$T_l(n,r)=\f{l!p^l}{\lfloor n/p^{-1}\rfloor!}\sum_{k=0}^n\bi nk(-1)^k\bi{k-r}l
=\f{l!p^l}{(pn)!}(-1)^n\bi{-r}{l-n},$$
where we have applied a known identity (cf. [GKP, (5.24)]) to get the last equality.
If $l<n$, then $T_l(n,r)=0\in\Z_p$. When $l\gs n$, we also have
$T_l(n,r)\in\Z_p$ because
$\ord_p((pn)!)=n+\ord_p(n!)\ls l+\ord_p(l!)$.
This ends the proof. \qed
\smallskip

For convenience, in the following lemmas
we let $p$ be a fixed prime and $\al$ be a positive integer.

\proclaim{Lemma 2.2} Let $l\in\N$, $n\in\Z^+$ and $r\in\Z$. Then
$$T_l(n-1,r)-T_l(n-1,r-1)=\cases T_l(n,r)&\t{if}\ p^{\al-1}\nmid n,
\\\f n{p^{\al-1}}T_l(n,r)&\t{otherwise}.\endcases\tag2.2$$
When $l>0$, we also have
$$\aligned&T_l(n,r)+\f r{p^{\al-1}}T_{l-1}(n,r+p^{\al})
\\=&\cases -T_{l-1}(n-1,r+p^{\al}-1)&\t{if}\ p^{\al-1}\mid n,
\\-\f n{p^{\al-1}}T_{l-1}(n-1,r+p^{\al}-1)&\t{otherwise}.
\endcases
\endaligned\tag2.3$$
\endproclaim
\Proof. Clearly
$$\align&\sum_{k\eq r\,(\mo\ p^{\al})}\bi nk(-1)^k\bi{(k-r)/p^{\al}}{l}
\\&+\sum_{k\eq r-1\,(\mo\ p^{\al})}\bi{n-1}k(-1)^k\bi{(k-(r-1))/p^{\al}}{l}
\\=&\sum_{k\eq r\,(\mo\ p^{\al})}\bi nk(-1)^k\bi{(k-r)/p^{\al}}{l}
\\&-\sum_{k\eq r\,(\mo\ p^{\al})}\bi{n-1}{k-1}(-1)^k\bi{(k-r)/p^{\al}}{l}
\\=&\sum_{k\eq r\,(\mo\ p^{\al})}\bi{n-1}k(-1)^k\bi{(k-r)/p^{\al}}{l}.
\endalign$$
Therefore
$$\l\lfloor\f n{p^{\al-1}}\r\rfloor!\,T_{l}(n,r)
+\l\lfloor\f{n-1}{p^{\al-1}}\r\rfloor!\,T_{l}(n-1,r-1)
=\l\lfloor\f{n-1}{p^{\al-1}}\r\rfloor!\,T_{l}(n-1,r)$$
and hence (2.2) follows.

Now let $l>0$.
Note that
$$\align&\l\lfloor \f n{p^{\al-1}}\r\rfloor!\f{T_l(n,r)}{(l-1)!p^{l-1}}
\\=&lp\sum_{k\eq r\,(\mo\ p^{\al})}\bi nk(-1)^k
\f {(k-r)/p^{\al}}l\bi{(k-r)/p^{\al}-1}{l-1}
\\=&\sum_{k\eq r\,(\mo\ p^{\al})}\bi nk(-1)^k\(\f{k}{p^{\al-1}}-\f r{p^{\al-1}}\)
\bi{(k-r-p^{\al})/p^{\al}}{l-1}.
\endalign$$
Using the identity $\bi nkk=n\bi{n-1}{k-1}$, we find that
$$\align&\l\lfloor \f n{p^{\al-1}}\r\rfloor!\f{T_l(n,r)}{(l-1)!p^{l-1}}
\\=&\f n{p^{\al-1}}\sum_{k\eq r\,(\mo\ p^{\al})}
\bi {n-1}{k-1}(-1)^k\bi{(k-r-p^{\al})/p^{\al}}{l-1}
\\&-\f r{p^{\al-1}}\sum_{k\eq r\,(\mo\ p^{\al})}
\bi {n}k(-1)^k\bi{(k-r-p^{\al})/p^{\al}}{l-1}
\\=&\f n{p^{\al-1}}\sum_{k\eq r-1\,(\mo\ p^{\al})}\bi {n-1}k(-1)^{k+1}
\bi{(k-(r-1)-p^{\al})/p^{\al}}{l-1}
\\&-\f r{p^{\al-1}}\sum_{k\eq r\,(\mo\ p^{\al})}\bi {n}k(-1)^k
\bi{(k-r-p^{\al})/p^{\al}}{l-1}.
\endalign$$
So we have
$$\align\l\lfloor \f n{p^{\al-1}}\r\rfloor!\, T_l(n,r)
=&-\f n{p^{\al-1}}\l\lfloor\f{n-1}{p^{\al-1}}\r\rfloor!
\,T_{l-1}(n-1,r+p^{\al}-1)
\\&-\f r{p^{\al-1}}\l\lfloor\f{n}{p^{\al-1}}\r\rfloor!
\,T_{l-1}(n,r+p^{\al}),
\endalign$$
which is equivalent to (2.3). \qed

\Remark\ 2.1. Lemma 2.2 is not sufficient for an induction proof of Theorem 2.1;
in fact we immediately encounter difficulty when $p^{\al-1}\mid n$ and $p^{\al-1}\nmid r$.

\proclaim{Lemma 2.3} Let $d,m\in\Z^+$, $n\in\N$ and $r\in\Z$, and let $f(x)$
be a function from $\Z$ to the complex field. Then we have
$$\aligned&\sum_{k\eq r\,(\mo\ d)}\bi nk(-1)^kf\l(\l\lfloor\f{k-r}m\r\rfloor\r)
\\=&\sum_{j=0}^{n}\bi nj\(\sum_{i\eq r\,(\mo\ d)}\bi ji(-1)^i\)
\sum_{i=0}^{m-1} \sigma_{ij},
\endaligned\tag2.4$$
where
$$\sigma_{ij}=\sum_{k\eq r+i-j\,(\mo\ m)}\bi{n-j}k(-1)^kf\l(\f{k-(r+i-j)}m\r).\tag2.5$$
\endproclaim

\Proof. Let $\zeta$ be a primitive $d$th root of unity.
Given $j\in\{0,\ldots,n\}$, we have
$$\align\sum_{i\eq r\,(\mo\ d)}\bi ji(-1)^i
=&\sum_{i=0}^j\bi ji\f{(-1)^i}d\sum_{s=0}^{d-1}\zeta^{(i-r)s}
\\=&\f1d\sum_{s=0}^{d-1}\zeta^{-rs}\sum_{i=0}^j\bi ji(-\zeta^s)^i
\\=&\f1d\sum_{s=0}^{d-1}\zeta^{-rs}(1-\zeta^s)^j.
\endalign$$
Also,
$$\align\sum_{i=0}^{m-1}\sigma_{ij}
=&\sum_{i=0}^{m-1}\sum_{k-(r-j)\eq i\,(\mo\ m)}\bi{n-j}k(-1)^k
f\l(\f{k-(r-j)-i}m\r)
\\=&\sum_{i=0}^{m-1}\sum_{k-(r-j)\eq i\,(\mo\ m)}\bi{n-j}k(-1)^k
f\l(\l\lfloor\f{k-(r-j)}m\r\rfloor\r)
\\=&\sum_{k=0}^{n-j}\bi{n-j}k(-1)^kf\l(\l\lfloor\f{k+j-r}m\r\rfloor\r).
\endalign$$
Therefore
$$\align&\sum_{j=0}^n\bi nj\(\sum_{i\eq r\,(\mo\ d)}\bi ji(-1)^i\)\sum_{i=0}^{m-1}\sigma_{ij}
\\=&\sum_{j=0}^n\bi nj\f1d\sum_{s=0}^{d-1}\zeta^{-rs}(1-\zeta^s)^j
\sum_{k=j}^n\bi{n-j}{k-j}(-1)^{k-j}f\l(\l\lfloor\f{k-r}m\r\rfloor\r)
\\=&\f1d\sum_{s=0}^{d-1}\zeta^{-rs}\sum_{k=0}^n\bi nk(-1)^kf\l(\l\lfloor\f{k-r}m\r\rfloor\r)
\sum_{j=0}^k\bi kj(\zeta^s-1)^j
\\=&\sum_{k=0}^n\bi nk(-1)^kf\l(\l\lfloor\f{k-r}m\r\rfloor\r)
\f1d\sum_{s=0}^{d-1}\zeta^{(k-r)s}
\\=&\sum_{k\eq r\,(\mo\ d)}\bi nk(-1)^kf\l(\l\lfloor\f{k-r}m\r\rfloor\r).
\endalign$$
This proves (2.4). \qed
\smallskip

\proclaim{Lemma 2.4} Let $l,n\in\N$ and $r\in\Z$. Then
$$T_l(n,r)=\sum_{i=0}^{p^{\al}-1}\sum_{j=0}^n c_{\al}(n,j)T_0(j,r)T_l(n-j,\,r+i-j),\tag2.6$$
where
$$c_{\al}(n,j):=\bi nj\f{\lfloor j/p^{\al-1}\rfloor!\lfloor(n-j)/p^{\al-1}\rfloor!}
{\lfloor n/p^{\al-1}\rfloor!}.\tag2.7$$
\endproclaim
\Proof. It suffices to apply Lemma 2.3 with $d=m=p^{\al}$ and $f(x)=\bi xl$. \qed

\proclaim{Lemma 2.5} We have $c_{\al}(n,j)\in\Z_p$ for all $j=0,\ldots,n$.
\endproclaim
\Proof. Clearly
$$\align\ord_p(c_{\al}(n,j))=&\ord_p(n!)-\ord_p\l(\l\lfloor\f n{p^{\al-1}}\r\rfloor!\r)
\\&-\(\ord_p(j!)-\ord_p\l(\l\lfloor\f j{p^{\al-1}}\r\rfloor!\r)\)
\\&-\(\ord_p((n-j)!)-\ord_p\l(\l\lfloor\f {n-j}{p^{\al-1}}\r\rfloor!\r)\)
\\=&\sum_{0<s<\al}\(\l\lfloor\f n{p^s}\r\rfloor-\l\lfloor\f j{p^s}\r\rfloor
-\l\lfloor\f{n-j}{p^s}\r\rfloor\)\gs0.
\endalign$$
This concludes the proof. \qed

\medskip
\noindent{\it Proof of Theorem 2.1}. Lemma 2.1 indicates that
Theorem 2.1 holds when $\al=0$. Below we let $\al>0$.

{\bf Step I}. Use induction on $l+n$ to show that $T_l(n,r)\in\Z_p$ for any $r\in\Z$.

The case $l=n=0$ is trivial since $T_0(0,r)\in\Z$.

Now let $l+n>0$, and assume that $T_{l_*}(n_*,r_*)\in\Z_p$
whenever $l_*,n_*\in\N$, $l_*+n_*<l+n$ and $r_*\in\Z$.

{\tt Case 1}. $l=0$. By Weisman's result mentioned in the first
section (see also, [S06]),
$$\align&\ord_p\(\sum_{k\eq r\,(\mo\ p^{\al})}\bi{n}k(-1)^k\)
\\\gs&\l\lfloor\f{n-p^{\al-1}}{\varphi(p^{\al})}\r\rfloor
=\l\lfloor\f{n/p^{\al-1}-1}{p-1}\r\rfloor=\l\lfloor\f{n_0-1}{p-1}\r\rfloor
\endalign$$
where $n_0=\lfloor n/p^{\al-1}\rfloor$. If $n_0>0$, then
$\ord_p(n_0!)\ls\lfloor(n_0-1)/(p-1)\rfloor$ (as mentioned in the first section);
hence
$$T_0(n,r)=\f1{n_0!}\sum_{k\eq r\,(\mo\ p^{\al})}\bi{n}k(-1)^k\in\Z_p.$$
Clearly this also holds when $n_0=0$.

{\tt Case 2}. $l>0$ and $p^{\al}\nmid r$.
In this case, $T_0(0,r)$ vanishes.
Thus, by Lemmas 2.4--2.5 and
the induction hypothesis, we have
$$T_l(n,r)=\sum_{i=0}^{p^{\al}-1}\sum_{0<j\ls n}c_{\al}(n,j)
T_0(j,r)T_l(n-j,\,r+i-j)\in\Z_p.$$

{\tt Case 3}. $l>0$ and $p^{\al}\mid r$.
If $p^{\al-1}\nmid n$, then
$$T_l(n,r)=T_l(n-1,r)-T_l(n-1,r-1)\in\Z_p$$
by (2.2) and the induction hypothesis.
If $p^{\al-1}\mid n$ and $n\not=0$, then
$$T_l(n,r)=-\f r{p^{\al-1}}T_{l-1}(n,r+p^{\al})-T_{l-1}(n-1,r+p^{\al}-1)\in\Z_p$$
by (2.3) and the induction hypothesis. Note also that $T_l(0,r)\in\Z_p$.

In view of the above, we have finished the first step.
\medskip

{\bf Step II}. Use induction on $n$ to prove (2.1) for any $r\in\Z$.

If $p^{\al}\nmid r$, then $T_l(0,r)=0$; if $p^{\al}\mid r$ then
$\tau_p(\{r\}_{p^{\al-1}},\{-r\}_{p^{\al-1}})=0$. So (2.1) holds when $n=0$.

Now let $n>0$ and $\tau_p(\{r\}_{p^{\al-1}},\{n-r\}_{p^{\al-1}})\not=0$.
Then both $\ord_p(r)$ and $\ord_p(n-r)$ are smaller than $\al-1$.
Assume that (2.1) with $n$ replaced by $n-1$ holds for all $r\in\Z$.
For $r'=n-r-(l-1)p^{\al}$, we have
$$\align\tau_p\l(\{r'\}_{p^{\al-1}},\{n-r'\}_{p^{\al-1}}\r)
=&\tau_p\l(\{n-r\}_{p^{\al-1}},\{r\}_{p^{\al-1}}\r)
\\=&\tau_p\l(\{r\}_{p^{\al-1}},\{n-r\}_{p^{\al-1}}\r)
\endalign$$
and
$$\align T_l(n,r')=&\f{l!p^l}{\lfloor n/p^{\al-1}\rfloor!}\sum_{k\eq r'\,(\mo\ p^{\al})}
\bi nk(-1)^k\bi{(n-r'-(n-k))/p^{\al}}l
\\=&\f{l!p^l}{\lfloor n/p^{\al-1}\rfloor!}\sum_{k\eq n-r'\,(\mo\ p^{\al})}
\bi nk(-1)^{n-k}\bi{(n-r'-k)/p^{\al}}l
\\=&\f{(-1)^{l+n} l!p^l}{\lfloor n/p^{\al-1}\rfloor!}\sum_{k\eq r\,(\mo\ p^{\al})}
\bi nk(-1)^{k}\bi{(k-(n-r'))/p^{\al}+l-1}l
\\=&(-1)^{l+n}T_l(n,r);
\endalign$$
also, $\ord_p(r')=\ord_p(n-r)=\ord_p(n)$ if $\ord_p(r)>\ord_p(n)$.
Thus, without loss of generality, below we simply let $\ord_p(r)\ls\ord_p(n)$.

In view of Lemma 2.2,
$$T_{l+1}(n,r-p^{\al})+\f{r-p^{\al}}{p^{\al-1}}T_l(n,r)
=\cases -T_l(n-1,r-1)&\t{if}\ p^{\al-1}\mid n,
\\-\f n{p^{\al-1}}T_l(n-1,r-1)&\t{otherwise}.\endcases$$
As $T_{l+1}(n,r-p^{\al})\in\Z_p$,
 and $$\ord_p(T_l(n-1,r-1))\gs\tau_p(\{r-1\}_{p^{\al-1}},\{n-1-(r-1)\}_{p^{\al-1}})$$
 by the induction hypothesis, we have
$$\align\ord_p(rT_l(n,r))&\gs\min\{\al-1,\ord_p(nT_l(n-1,r-1))\}
\\&\gs\min\{\al-1,\ord_p(n)+\tau_p(\{r-1\}_{p^{\al-1}},\{n-r\}_{p^{\al-1}})\}.
\endalign$$
By the definition of $\tau_p$ and the inequality $\ord_p(r)\ls\ord_p(n)$,
$$\tau_p(\{r\}_{p^{\al-1}},\{n-r\}_{p^{\al-1}})\ls\al-1-\ord_p(r)$$
and also
$$\tau_p(\{r\}_{p^{\al-1}},\{n-r\}_{p^{\al-1}})
\ls\tau_p(\{r-1\}_{p^{\al-1}},\{n-r\}_{p^{\al-1}})+\ord_p(n)-\ord_p(r).$$
(Note that $\{r\}_{p^{\al-1}}+\{n-r\}_{p^{\al-1}}\eq n\ (\mo\ p^{\al-1})$.)
So, (2.1) follows from the above.

The proof of Theorem 2.1 is now complete. \qed

\heading{3. Proofs of Theorem 1.7 and Corollary 1.2}\endheading

To prove Theorem 1.7 and Corollary 1.2, we need to establish an auxiliary theorem first.

\proclaim{Lemma 3.1} Let $m\in\Z^+$ and $r\in\Z$ be relatively prime.
Then, for any $n\in\Z^+$, we have
$$\f1n\sum_{k=0}^{n-1}\f1{km+r}\eq\f1r+\f m2[\![2\mid n]\!]\ \ (\mo\ m).\tag3.1$$
\endproclaim
\Proof. We use induction on $n$.

If $n\in\Z^+$ is relatively prime to $m$ (e.g., $n=1$),
then $2$ cannot divide both $m$ and $n$; hence
$$\f1n\sum_{k=0}^{n-1}\f1{km+r}\eq\f1n\sum_{k=0}^{n-1}\f1r
\eq\f1r+\f m2[\![2\mid n]\!]\ \ (\mo\ m).$$

Now suppose that $p$ is the least common prime divisor of $m$ and $n$, and set $n_0=n/p$. Then
$$\sum_{k=0}^{n-1}\f1{km+r}=\sum_{i=0}^{n_0-1}\sum_{j=0}^{p-1}\f1{(i+jn_0)m+r}
=\sum_{i=0}^{n_0-1}\sum_{j=0}^{p-1}\f1{im+jmn_0+r}.$$
For any $i=0,\ldots,n_0-1$, clearly
$$\align&\f1n\sum_{j=0}^{p-1}\(\f1{im+jmn_0+r}-\f1{im+r}\)
\\=&\f1n\sum_{j=0}^{p-1}\f{-jmn/p}{(im+jmn_0+r)(im+r)}
\\\eq&-\sum_{j=0}^{p-1}\f {jm/p}{r^2}=-\f{p-1}{2}\cdot\f m{r^2}
\eq\da_{p,2}\f m2\ \ (\mo\ m).
\endalign$$
Therefore
$$\align\f1n\sum_{k=0}^{n-1}\f1{km+r}
\eq&\sum_{i=0}^{n_0-1}\(\f1n\sum_{j=0}^{p-1}\f1{im+r}+\da_{p,2}\f m2\)
\\\eq&\f1{n_0}\sum_{i=0}^{n_0-1}\f1{im+r}+\da_{p,2}\f m2\cdot\f n2\ \ (\mo\ m).
\endalign$$
Note that $n_0<n$. If
$$\f1{n_0}\sum_{i=0}^{n_0-1}\f1{im+r}\eq\f1r+\f m2[\![2\mid n_0]\!]\ \ (\mo\ m),$$
then (3.1) holds by the above, because
$$\da_{p,2}\f m2\cdot\f n2+\f m2[\![2\mid n_0]\!]
\eq\da_{p,2}\f m2([\![4\mid n-2]\!]+[\![4\mid n]\!])\eq\f m2[\![2\mid n]\!] \ (\mo\ m).$$
This concludes the induction proof. \qed

\proclaim{Lemma 3.2} Let $p$ be a prime, and let $k\in\N$ and $n\in\Z^+$.
If $p$ is odd, then
$$\bi{pn}{pk}\eq\bi nk\ \l(\mo\ p^{2\ord_p(n)+2}\r).\tag3.2$$
For $p=2$ we have
$$\bi{2n}{2k}\eq(-1)^k\bi nk\ \l(\mo\ 2^{2\ord_2(n)+1}\r).\tag3.3$$
\endproclaim
\Proof. The case $k=0$ or $k\gs n$ is trivial. Below we let $0<k<n$.

 By a result of Jacobsthal (see, e.g., [G]), if $p>3$, then
$$\bi{pn}{pk}\bigg/\bi nk=1+p^3nk(n-k)q$$
for some $q\in\Z_p$, and hence
$$\align \bi{pn}{pk}-\bi nk=&\bi nkp^3nk(n-k)q=p^3n^2\bi{n-1}{k-1}(n-k)q
\\\eq&0\ \l(\mo\ p^{3+2\ord_p(n)}\r).
\endalign$$

Now we handle the case $p=3$. Observe that
$$\align\bi{3n}{3k}\bigg/\bi nk
=&\f{(3n-1)(3n-2)}{1\cdot 2}\times\f{(3n-4)(3n-5)}{4\cdot 5}
\\&\times\cdots
\times\f{(3n-(3k-2))(3n-(3k-1))}{(3k-2)(3k-1)}
\endalign$$
and
$$(3n-(3i+1))(3n-(3i+2))=9n^2-9n(2i+1)+(3i+1)(3i+2)$$
for any $i\in\Z$. So we have
$$\align\bi{3n}{3k}\bigg/\bi nk=&\prod_{i=0}^{k-1}\(1+9n\f{n-2i-1}{(3i+1)(3i+2)}\)
\\\eq&1+9n\sum_{i=0}^{k-1}\f{n-2i-1}{(3i+1)(3i+2)}\ \l(\mo\ (3^{2+\ord_3(n)})^2\r).
\endalign$$
Clearly,
$$\align\sum_{i=0}^{k-1}\f{n-2i-1}{(3i+1)(3i+2)}=&\sum_{i=0}^{k-1}\(\f{n-2i-1}{3i+1}
-\f{n-2i-1}{3i+2}\)
\\=&\sum_{i=0}^{k-1}\(\f{n-2i-1}{3i+1}
-\f{n-2(k-1-i)-1}{3(k-1-i)+2}\)
\endalign$$
and hence
$$\align&\f1k\sum_{i=0}^{k-1}\f{n-2i-1}{(3i+1)(3i+2)}
\\=&\f1k\sum_{i=0}^{k-1}\(\f{n-(2i+1)}{3i+1}-\f{n+(2i+1)-2k}{3k-(3i+1)}\)
\\=&\f1k\sum_{i=0}^{k-1}\f{3k(n-(2i+1))+2k(3i+1)-2n(3i+1)}{(3i+1)(3k-3i-1)}
\\=&\sum_{i=0}^{k-1}\f{3n-1}{(3i+1)(3k-3i-1)}-\f{2n}k\sum_{i=0}^{k-1}\f1{3(k-1-i)+2}
\\\eq&\sum_{i=0}^{k-1}\f{-1}{-1}-2n\times\f12=k-n\ \ (\mo\ 3),
\endalign$$
where we have applied Lemma 3.1 with $m=3$ to get the last congruence.
Therefore
$$\align\bi{3n}{3k}\eq&\bi nk+9n\f nk\bi{n-1}{k-1}
\sum_{i=0}^{k-1}\f{n-2i-1}{(3i+1)(3i+2)}\ \l(\mo\ 3^{2\ord_3(n)+4}\r)
\\\eq&\bi nk+9n^2\bi{n-1}{k-1}(k-n)\ \l(\mo\ 3^{2\ord_3(n)+3}\r)
\\\eq&\bi nk\ \l(\mo\ 3^{2\ord_3(n)+2}\r).
\endalign$$

Finally we consider the case $p=2$.
Observe that
$$\align\f{\bi{2n}{2k}}{\bi{n}{k}}=&\prod_{j=0}^{k-1}\f{2n-(2j+1)}{2j+1}
=(-1)^k\prod_{j=0}^{k-1}\l(1-\f{2n}{2j+1}\r)
\\\eq&(-1)^k\(1-2n\sum_{j=0}^{k-1}\f1{2j+1}\)\ \l(\mo\ (2^{\ord_2(n)+1})^2\r).
\endalign$$
This, together with Lemma 3.1 in the case $m=2$, yields that
$$\align\bi{2n}{2k}\eq&(-1)^k\(\bi nk-\f{2n^2}k\bi{n-1}{k-1}\sum_{j=0}^{k-1}\f1{2j+1}\)
\ \l(\mo\ 2^{2\ord_2(n)+2}\r)
\\\eq&(-1)^k\bi nk-(-1)^k2n^2\bi{n-1}{k-1}(1+[\![2\mid k]\!])\ \ \l(\mo\ 2^{2\ord_2(n)+2}\r)
\\\eq&(-1)^k\bi nk\ \ \l(\mo\ 2^{2\ord_2(n)+1}\r).
\endalign$$
We are done. \qed
\medskip

If $p$ is a prime, $\al\in\Z^+$, $n\in\N$ and $r\in\Z$, then
$$\ord_p\(\sum_{k\eq r\,(\mo\ p^{\al})}\bi{p^{\al-1}n}k(-1)^k\)
\gs\l\lfloor\f{p^{\al-1}n-p^{\al-1}}{\varphi(p^{\al})}\r\rfloor
=\l\lfloor\f{n-1}{p-1}\r\rfloor$$
by Weisman's result, and hence
$$S_{\al}^{(p)}(n,r):=p^{-\lfloor\f{n-1}{p-1}\rfloor}
\sum_{k\eq r\,(\mo\ p^{\al})}\bi{p^{\al-1}n}k(-1)^k\tag3.4$$
is an integer.

\proclaim{Theorem 3.1} Let $p$ be a prime.
Then, for each $\al=2,3,\ldots$, whenever $n\in\N$ and $r\in\Z$ we have
$$S_{\al}^{(p)}(n,r)\eq\cases S_{\al-1}^{(p)}(n,r/p)
\ (\mo\ p^{(2-\da_{p,2})(\al-2)})&\t{if}\ p\mid r,
\\0\ (\mo\ p^{\al-2})\ &\t{otherwise},\endcases\tag3.5$$
where $S_{\al}^{(p)}(n,r)$ is defined by $(3.4)$.
\endproclaim
\Proof. We use induction on $\al$ and write $S_{\al}(n,r)$ for $S_{\al}^{(p)}(n,r)$.

When $\al=2$ we need do nothing.

Now let $\al>2$ and assume the desired result for $\al-1$.
We use induction on $n$ to prove that (3.5) holds for any $r\in\Z$.

First we consider the case $n\in\{0,\ldots,p-1\}$.
Note that $p^{\al-1}n<p^{\al}$ and $\da=-\lfloor(n-1)/(p-1)\rfloor\in\{0,1\}$.
Thus
$$S_{\al}(n,r)=p^{\da}\bi{p^{\al-1}n}{\{r\}_{p^{\al}}}(-1)^{\{r\}_{p^{\al}}}.$$
If $p\nmid r$, then
$$\bi{p^{\al-1}n}{\{r\}_{p^{\al}}}=\f{p^{\al-1}n}{\{r\}_{p^{\al}}}
\bi{p^{\al-1}n-1}{\{r\}_{p^{\al}}-1}\eq0\ (\mo\ p^{\al-1}),$$
and hence $S_{\al}(n,r)\eq0\ (\mo\ p^{\al-1})$.
When $p\mid r$, by Lemma 3.2 we have
$$\align&(-1)^{\{r\}_{p^{\al}}-\{r/p\}_{p^{\al-1}}}\bi{p^{\al-2}n}{\{r/p\}_{p^{\al-1}}}
\\=&(-1)^{(p-1)\{r\}_{p^{\al}}/p}\bi{p^{\al-2}n}{\{r\}_{p^{\al}}/p}
\eq\bi{p^{\al-1}n}{\{r\}_{p^{\al}}}
\ \l(\mo\ p^{2(\al-2)+2-\da_{p,2}}\r),
\endalign$$
thus $\ord_p(S_{\al}(n,r)-S_{\al-1}(n,r/p))\gs 2\al-2-\da_{p,2}>(2-\da_{p,2})(\al-2)$.

Below we let $n\gs p$, and assume that (3.5) is valid for all $r\in\Z$
if we replace $n$ in (3.5) by a smaller nonnegative
integer. Set $n'=n-(p-1)\gs1$. Then,
by Vandermonde's identity $\bi{x+y}k=\sum_{j\in\N}\bi xj\bi y{k-j}$
(cf. [GKP, (5.22)]),
we have
$$\align S_{\al}(n,r)=&p^{-\lfloor\f{n-1}{p-1}\rfloor}\sum_{k\eq r\,(\mo\ p^{\al})}
\sum_{j=0}^{\varphi(p^{\al})}\bi{\varphi(p^{\al})}j\bi{p^{\al-1}n'}{k-j}(-1)^k
\\=&\sum_{j=0}^{\varphi(p^{\al})}\f{(-1)^j}p\bi{\varphi(p^{\al})}j
p^{-\lfloor\f{n'-1}{p-1}\rfloor}\sum_{k\eq r-j\,(\mo\ p^{\al})}\bi{p^{\al-1}n'}k(-1)^k.
\endalign$$
If $0\ls j\ls\varphi(p^{\al})$ and $p\nmid j$, then
$$\f 1p\bi{\varphi(p^{\al})}j=\f{p^{\al-1}(p-1)}{pj}
\bi{\varphi(p^{\al})-1}{j-1}\eq0\ (\mo\ p^{\al-2}),$$
and also $S_{\al}(n',r-j)\eq0\ (\mo\ p^{\al-2})$
providing $p\mid r$ (by the induction hypothesis). Thus
$$S_{\al}(n,r)\eq\sum_{j=0}^{\varphi(p^{\al-1})}\f{(-1)^{pj}}p
\bi{\varphi(p^{\al})}{pj}S_{\al}(n',r-pj)\ \l(\mo\ p^{(1+[\![p\mid r]\!])(\al-2)}\r).$$
Note that when
$j\not\eq0\ (\mo\ p^{\al-2})$ we have
$$\f1p\bi{\varphi(p^{\al})}{pj}=\f{p^{\al-1}(p-1)}{p^2j}\bi{\varphi(p^{\al})-1}{pj-1}
\in\Z_p.$$

{\tt Case I}. $p\nmid r$. By the induction hypothesis,
 $p^{\al-2}\mid S_{\al}(n',r-pj)$ for all $j\in\Z$. Thus, by the above,
$$S_{\al}(n,r)\eq\sum_{j=0}^{p-1}\f{(-1)^{p^{\al-1}j}}p\bi{\varphi(p^{\al})}{p^{\al-1}j}
S_{\al}(n',r-p^{\al-1}j)\ (\mo\ p^{\al-2}).$$
In view of Lucas' theorem or Lemma 3.2,
$$\bi{p^{\al-1}(p-1)}{p^{\al-1}j}\eq\bi{p^{\al-2}(p-1)}{p^{\al-2}j}\eq\cdots
\eq\bi {p-1}j\eq(-1)^j\ (\mo\ p)$$
for every $j=0,\ldots,p-1$.
Note also that $(-1)^j\eq 1\ (\mo\ 2)$. So we have
$$S_{\al}(n,r)\eq\f1p\sum_{j=0}^{p-1}
S_{\al}(n',r-p^{\al-1}j)\ \ (\mo\ p^{\al-2}).$$
Observe that
$$\align&\sum_{j=0}^{p-1}S_{\al}(n',r-p^{\al-1}j)
\\=&p^{-\lfloor\f{n'-1}{p-1}\rfloor}\sum_{k\eq r\,(\mo\ p^{\al-1})}\bi{p^{\al-2}pn'}k(-1)^k
\\=&p^{-\lfloor\f{n'-1}{p-1}\rfloor}p^{\lfloor\f{pn'-1}{p-1}\rfloor}S_{\al-1}(pn',r)
=p^{n'}S_{\al-1}(pn',r).
\endalign$$
Therefore
$$S_{\al}(n,r)\eq p^{n'-1}S_{\al-1}(pn',r)\ \l(\mo\ p^{\al-2}\r).$$

 By the induction hypothesis for $\al-1$, we have $p^{\al-3}\mid S_{\al-1}(pn',r)$, hence
$S_{\al}(n,r)\eq0\ (\mo\ p^{\al-2})$ if $n'>1$. In the case $n'=1$, we need to show that
$S_{\al-1}(p,r)\eq0\ (\mo\ p^{\al-2})$. In fact,
$$\align S_{\al-1}(p,r)=&p^{-\lfloor\f{p-1}{p-1}\rfloor}\sum_{k\eq r\,(\mo\ p^{\al-1})}
\bi{p^{\al-2}p}k(-1)^k
\\=&\sum_{k\eq r\,(\mo\ p^{\al-1})}\f{p^{\al-2}}k\bi{p^{\al-1}-1}{k-1}(-1)^k
\eq0\ (\mo\ p^{\al-2}).\endalign$$
(Note that if $k\eq r\ (\mo\ p^{\al-1})$ then $p\nmid k$ since $p\nmid r$.)

{\tt Case II}. $p\mid r$. By the paragraph before Case I,
$$\align &S_{\al}(n,r)-S_{\al-1}\l(n,\f rp\r)
\\\eq&\sum_{j=0}^{\varphi(p^{\al-1})}\f{(-1)^{pj}}p
\bi{\varphi(p^{\al})}{pj}S_{\al}(n',r-pj)
\\&-\sum_{j=0}^{\varphi(p^{\al-1})}\f{(-1)^{j}}p
\bi{\varphi(p^{\al-1})}{j}S_{\al-1}\l(n',\f rp-j\r)
\ \l(\mo\ p^{2(\al-2)}\r).
\endalign$$
By Lemma 3.2,
$$\bi{p^{\al-1}(p-1)}{pj}\eq(-1)^{(p-1)j}
\bi{p^{\al-2}(p-1)}{j}\ \l(\mo\ p^{2(\al-2)+1}\r).$$
Thus $S_{\al}(n,r)-S_{\al-1}(n,r/p)$ is congruent to
$$\sum_{j=0}^{\varphi(p^{\al-1})}\f{(-1)^{j}}p
\bi{\varphi(p^{\al-1})}{j}\l(S_{\al}(n',r-pj)-S_{\al-1}\l(n',\f rp-j\r)\r)$$
modulo $p^{2(\al-2)}$. If $p^{\al-2}\nmid j$, then
$$\f 1p\bi{\varphi(p^{\al-1})}j=\f{p^{\al-2}(p-1)}{pj}\bi{\varphi(p^{\al-1})-1}{j-1}\in\Z_p.$$
Since $S_{\al}(n',r-pj)\eq S_{\al-1}(n',r/p-j)\ (\mo\ p^{(2-\da_{p,2})(\al-2)})$
for any $j\in\Z$ (by the induction hypothesis), and $-1\eq1\ (\mo\ 2)$,
we have
$$\align&S_{\al}(n,r)-S_{\al-1}\l(n,\f rp\r)
\\\eq&\sum_{j=0}^{p-1}\f{(-1)^{j}}p
\bi{\varphi(p^{\al-1})}{p^{\al-2}j}
\l(S_{\al}(n',r-p^{\al-1}j)-S_{\al-1}\l(n',\f rp-p^{\al-2}j\r)\r)
\\\eq&\sum_{j=0}^{p-1}\f{S_{\al}(n',r-p^{\al-1}j)-S_{\al-1}(n',r/p-p^{\al-2}j)}p
\ \l(\mo\ p^{(2-\da_{p,2})(\al-2)}\r),
\endalign$$
where we have used the fact that
$\bi{p^{\al-2}(p-1)}{p^{\al-2}j}\eq\bi{p-1}j\eq(-1)^j\ (\mo\ p)$
for $j=0,\ldots,p-1$.
As in Case I,
$$\sum_{j=0}^{p-1}S_{\al}(n',r-p^{\al-1}j)=p^{n'}S_{\al-1}(pn',r)$$
and
$$\sum_{j=0}^{p-1}S_{\al-1}\l(n',\f rp-p^{\al-2}j\r)
=p^{n'}S_{\al-2}\l(pn',\f rp\r).$$
Thus
$$\align& S_{\al}(n,r)-S_{\al-1}\l(n,\f rp\r)
\\\eq&p^{n'-1}\l(S_{\al-1}(pn',r)-S_{\al-2}\l(pn',\f rp\r)\r)
\ \l(\mo\ p^{(2-\da_{p,2})(\al-2)}\r).
\endalign$$

By the induction hypothesis for $\al-1$,
$S_{\al-1}(pn',r)$ and $S_{\al-2}(pn',r/p)$ are congruent
modulo $p^{(2-\da_{p,2})(\al-3)}$.
Therefore, if $n'>2$ then
$$S_{\al}(n,r)\eq S_{\al-1}\l(n,\f rp\r)\ \l(\mo\ p^{(2-\da_{p,2})(\al-2)}\r).$$
When $n'\in\{1,2\}$, we have $n'-1-\lfloor(pn'-1)/(p-1)\rfloor\gs-2$ and also
$$\align
&p^{\lfloor\f{pn'-1}{p-1}\rfloor}\l(S_{\al-1}(pn',r)-S_{\al-2}\l(pn',\f rp\r)\r)
\\=&\sum_{k\eq r\,(\mo\ p^{\al-1})}\bi{p^{\al-1}n'}k(-1)^k
-\sum_{k\eq r/p\,(\mo\ p^{\al-2})}\bi{p^{\al-2}n'}k(-1)^k
\\=&\sum_{k\eq r/p\,(\mo\ p^{\al-2})}\(\bi{p^{\al-1}n'}{pk}(-1)^{pk}
-\bi{p^{\al-2}n'}{k}(-1)^k\).
\endalign$$
By Lemma 3.2, the last sum is a multiple of
$$p^{2(\al-2)+2-\da_{p,2}}=p^{(2-\da_{p,2})(\al-2)}p^{2+\da_{p,2}(\al-3)}.$$
So we also have the desired result in the case $n'\ls 2$.

The induction proof of Theorem 3.1 is now complete. \qed
\smallskip

We believe that (3.5) in the case $p\mid r$ can be improved slightly. Here is our conjecture.
\proclaim{Conjecture 3.1} Let $p$ be an odd prime, and let $\al\gs2$ be an integer. If
$n\in\N$ and $r\in\Z$, then
$$S_{\al}^{(p)}(n,pr)\eq S_{\al-1}^{(p)}(n,r)\ \ \l(\mo\ p^{2\al-2-\da_{p,3}}\r).$$
\endproclaim

Now we give a useful consequence of Theorem 3.1.
\proclaim{Corollary 3.1} Let $p$ be a prime, and let $\al,\beta,n\in\N$ with $\al>\beta$.
Given $r\in\Z$ we have
$$S_{\al}^{(p)}(n,p^{\beta}r)\eq S_{\al-\beta}^{(p)}(n,r)
\ \ \l(\mo\ p^{(2-\da_{p,2})(\al-\beta-1)}\r).\tag3.6$$
Provided that $r\not\eq0\ (\mo\ p)$, we also have
$$\ord_p\l(S_{\al}^{(p)}(n,p^{\beta}r)\r)\gs\al-\beta-2.\tag3.7$$
\endproclaim
\Proof. By Theorem 3.1, if $0\ls j<\beta$ then
$$S_{\al-j}^{(p)}(n,p^{\beta-j}r)\eq S_{\al-j-1}^{(p)}(n,p^{\beta-j-1}r)
\ \l(\mo\ p^{(2-\da_{p,2})(\al-j-2)}\r).$$
So (3.6) follows.

When $\al=\beta+1$, (3.7) holds trivially.
 If $p\nmid r$ and $\al-\beta\gs2$,
then $S_{\al-\beta}^{(p)}(n,r)\eq0\ (\mo\ p^{\al-\beta-2})$ by Theorem 3.1,
 combining this with (3.6) we immediately obtain (3.7). \qed

\medskip
\noindent {\it Proof of Theorem 1.7}. Write $n=pn_0+s_0$ with $n_0,s_0\in\N$ and $s_0<p$.
Then
$$\align\l\lfloor\f{p^{\al-2}n+s-p^{\al-1}}
{\varphi(p^{\al})}\r\rfloor
=&\l\lfloor \f{n+s/p^{\al-2}-p}{\varphi(p^{2})}\r\rfloor
=\l\lfloor \f{n-p}{\varphi(p^{2})}\r\rfloor
\\=&\l\lfloor\f{n/p-1}{p-1}\r\rfloor
=\l\lfloor\f{n_0-1}{p-1}\r\rfloor.
\endalign$$
By Vandermonde's identity,
$$\align&\sum_{k\eq p^{\al-2}r+t\,(\mo\ p^{\al})}\bi{p^{\al-2}n+s}k(-1)^k
\\=&\sum_{k\eq p^{\al-2}r+t\,(\mo\ p^{\al})}\sum_{j\in\N}\bi{p^{\al-2}s_0+s}j
\bi{p^{\al-1}n_0}{k-j}(-1)^k
\\=&\sum_{j\in\N}\bi{p^{\al-2}s_0+s}j(-1)^j
\sum_{k\eq p^{\al-2}r+t-j\,(\mo\ p^{\al})}\bi{p^{\al-1}n_0}{k}(-1)^k.
\endalign$$

In light of Corollary 3.1, for any $r'\in\Z$ we have
$$S_{\al}(n_0,r')\eq\cases S_2(n_0,r'/p^{\al-2})\ (\mo\ p^{2-\da_{p,2}})
&\t{if}\ r'\eq0\ (\mo\ p^{\al-2}),\\0\ (\mo\ p)&\t{if}\ \ord_p(r')<\al-2.\endcases$$
Thus
$$\align&p^{-\lfloor\f{n_0-1}{p-1}\rfloor}
\sum_{k\eq p^{\al-2}r+t\,(\mo\ p^{\al})}\bi{p^{\al-2}n+s}k(-1)^k
\\\eq&\sum_{i\in\N}\bi{p^{\al-2}s_0+s}{p^{\al-2}i+t}(-1)^{p^{\al-2}i+t}
p^{-\lfloor\f{n_0-1}{p-1}\rfloor}\sum_{k\eq r-i\,(\mo\ p^2)}\bi{pn_0}k(-1)^k
\\\eq&\sum_{i\in\N}\bi{s_0}i\bi st(-1)^{i+t}
p^{-\lfloor\f{n_0-1}{p-1}\rfloor}\sum_{k\eq r-i\,(\mo\ p^2)}\bi{pn_0}k(-1)^k
\ (\mo\ p),
\endalign$$
where we have applied Lucas' theorem to get the last congruence.
(Note that $s=t=0$ if $\al=2$, also $-1\eq 1\ (\mo\ 2)$.)
Since
$$\align&\sum_{i=0}^{s_0}\bi {s_0}i(-1)^i\sum_{k\eq r-i\,(\mo\ p^2)}\bi{pn_0}k(-1)^k
\\=&\sum_{i=0}^{s_0}\bi {s_0}i\sum_{k\eq r\,(\mo\ p^2)}\bi{pn_0}{k-i}(-1)^k
\\=&\sum_{k\eq r\,(\mo\ p^2)}\bi{pn_0+s_0}k(-1)^k=\sum_{k\eq r\,(\mo\ p^2)}\bi{n}k(-1)^k,
\endalign$$
the desired result follows from the above. \qed

\medskip
\noindent{\it Proof of Corollary 1.2}. Let $s=\{n\}_{2^{\al-2}}$ and $t=\{r\}_{2^{\al-2}}$.
Then $n=2^{\al-2}n_*+s$ and $r=2^{\al-2}r_*+t$.
By Theorem 1.7,
$$\align&2^{-\l\lfloor\f{n-2^{\al-1}}{\varphi(2^{\al})}\r\rfloor}
\sum_{k\eq r\,(\mo\ 2^{\al})}\bi nk(-1)^k\eq1\ (\mo\ 2)
\\\iff&\bi st\eq 2^{-\l\lfloor\f{n_*-2}{\varphi(4)}\r\rfloor}
\sum_{k\eq r_*\,(\mo\ 4)}\bi{n_*}k(-1)^k\eq1\ (\mo\ 2).
\endalign$$

In the case $n_*\eq1\ (\mo\ 2)$, by [S02, (3.3)] we have
$$2\sum\Sb 0<k<n_*\\4\mid k-r_*\endSb\bi {n_*}k
-\l(2^{n_*-1}-1\r)
=(-1)^{\f{r_*(n_*-r_*)}2}\l((-1)^{\f{n_*^2-1}8}2^{\f{n_*-1}2}-1\r),
$$
thus
$$\align2\sum_{k\eq r_*\,(\mo\ 4)}\bi {n_*}k
=&2([\![4\mid r_*]\!]+[\![4\mid n_*-r_*]\!])+2^{n_*-1}-1
\\&+(-1)^{\f{r_*(n_*-r_*)}2}\l((-1)^{\f{n_*^2-1}8}2^{\f{n_*-1}2}-1\r)
\\=&2^{n_*-1}+(-1)^{\f{r_*(n_*-r_*)}2+\f{n_*^2-1}8}2^{\f{n_*-1}2}
\endalign$$
and hence
$$2^{-\lfloor\f {n_*}2\rfloor+1}\sum_{k\eq r_*\,(\mo\ 4)}\bi {n_*}k\eq1\ (\mo\ 2)
\iff n_*>1.$$
When $n_*\eq0\ (\mo\ 2)$, if $n_*>0$ then by the above we have
$$\align&2\sum_{k\eq r_*\,(\mo\ 4)}\bi {n_*}k
=2\sum_{k\eq r_*\,(\mo\ 4)}\(\bi {n_*-1}k+\bi {n_*-1}{k-1}\)
\\=&2\sum_{k\eq r_*\,(\mo\ 4)}\bi {n_*-1}k+2\sum_{k\eq r_*-1\,(\mo\ 4)}\bi {n_*-1}k
\\=&2^{(n_*-1)-1}+(-1)^{\f{r_*(n_*-1-r_*)}2+\f{(n_*-1)^2-1}8}2^{\f{(n_*-1)-1}2}
\\&+2^{(n_*-1)-1}+(-1)^{\f{(r_*-1)(n_*-1-(r_*-1))}2+\f{(n_*-1)^2-1}8}2^{\f{(n_*-1)-1}2}
\\=&2^{n_*-1}+(-1)^{\bi{n_*/2}2}2^{n_*/2-1}
\l((-1)^{\f{r_*(n_*-1-r_*)}2}+(-1)^{\f{(r_*-1)(n_*-r_*)}2}\r)
\\=&2^{n_*-1}+(-1)^{\bi{n_*/2}2+\lfloor\f {r_*}2\rfloor}
\l(1+(-1)^{\f{n_*}2+r_*}\r)2^{n_*/2-1};
\endalign$$
therefore
$$\align&2^{-\lfloor\f {n_*}2\rfloor+1}\sum_{k\eq r_*\,(\mo\ 4)}\bi {n_*}k\eq1\ (\mo\ 2)
\\\iff& n_*>2\ \&\ \f{n_*}2\eq r_*\ (\mo\ 2),\ \t{or}\ n_*=2\ \&\ 2\mid r_*.
\endalign$$

Combining the above we obtain the desired result. \qed

\heading{4. Proofs of Theorems 1.5 and 1.8}\endheading

\noindent {\it Proof of Theorem 1.5}.
For convenience, we set $T_{l,\al}(n,r):=T_{l,\al}^{(p)}(n,r)$.

The case $n=0$ is easy, because
$$\align T_{l,\al+1}(0,r)=&[\![p^{\al+1}\mid r]\!]l!p^l\bi{-r/p^{\al+1}}l
\\=&[\![p\mid r\ \&\ p^{\al}\mid\lfloor r/p\rfloor]\!]l!p^l
\bi{-\lfloor r/p\rfloor/p^{\al}}l
\\=&(-1)^{\{r\}_p}\bi{0}{\{r\}_p}
T_{l,\al}\l(0,\l\lfloor \f rp\r\rfloor\r).
\endalign$$

Below we use induction on $l+n$ to prove the desired result.

If $l+n=0$, then $n=0$ and hence we are done.

Now let $n>0$, and assume that
$$T_{l_*,\al+1}(n_*,r_*)\eq(-1)^{\{r_*\}_p}\bi{\{n_*\}_p}{\{r_*\}_p}
T_{l_*,\al}\l(\l\lfloor \f {n_*}p\r\rfloor,
\l\lfloor\f {r_*}p\r\rfloor\r)\ (\mo\ p)$$
whenever $l_*,n_*\in\N$, $l_*+n_*<l+n$ and $r_*\in\Z$.
Write $n=pn_0+s$ and $r=pr_0+t$, where $n_0,r_0\in\Z$ and $s,t\in\{0,\ldots,p-1\}$.

{\tt Case 1}. $p^{\al}\nmid n$.
By Lemma 2.2 and the induction hypothesis, if $s\not=0$ then
$$\align T_{l,\al+1}(n,r)=&T_{l,\al+1}(pn_0+s-1,pr_0+t)
\\&-T_{l,\al+1}(pn_0+s-1,pr_0+t-1)
\\\eq&(-1)^t\bi{s-1}{t}T_{l,\al}(n_0,r_0)
\\&-\cases(-1)^{t-1}\bi{s-1}{t-1}T_{l,\al}(n_0,r_0)&\t{if}\ t>0,
\\(-1)^{p-1}\bi {s-1}{p-1}T_{l,\al}(n_0,r_0-1)&\t{if}\ t=0,\endcases
\\\eq&(-1)^t\bi{s-1}{t}T_{l,\al}(n_0,r_0)+(-1)^t\bi{s-1}{t-1}T_{l,\al}(n_0,r_0)
\\\eq&(-1)^t\bi stT_{l,\al}(n_0,r_0)\ \ (\mo\ p).
\endalign$$
Now we assume that $s=0$. In a similar way,
$$\align T_{l,\al+1}(n,r)=&T_{l,\al+1}(p(n_0-1)+p-1,pr_0+t)
\\&-T_{l,\al+1}(p(n_0-1)+p-1,pr_0+t-1)
\\\eq&(-1)^t\bi{p-1}{t}T_{l,\al}(n_0-1,r_0)
\\&-\cases(-1)^{t-1}\bi{p-1}{t-1}T_{l,\al}(n_0-1,r_0)\ (\mo\ p)&\t{if}\ t>0,
\\(-1)^{p-1}\bi {p-1}{p-1}T_{l,\al}(n_0-1,r_0-1)\ (\mo\ p)&\t{if}\ t=0.\endcases
\endalign$$
Thus, if $0<t<p$ then
$$T_{l,\al+1}(n,r)\eq(-1)^t\bi{p}{t}T_{l,\al}(n_0-1,r_0)
\eq0\ \  (\mo\ p)$$
and hence $T_{l,\al+1}(n,r)\eq(-1)^t\bi{0}{t}T_{l,\al}(n_0,r_0)\ (\mo\ p)$;
if $t=0$ then
$$T_{l,\al+1}(n,r)
\eq T_{l,\al}(n_0-1,r_0)-T_{l,\al}(n_0-1,r_0-1)=T_{l,\al}(n_0,r_0)\ (\mo\ p)$$
with the help of Lemma 2.2. (Note that $p^{\al-1}$ does not divide $n_0=n/p$.)

{\tt Case 2}. $p^{\al}\mid n$ and $p^{\al}\nmid r$. In this case,
$$\ord_p(T_{l,\al+1}(n,r))\gs\tau_p(\{r\}_{p^{\al}},\{n-r\}_{p^{\al}})>0$$
by Theorem 2.1.
If $p\nmid r$ (i.e., $t>0$) then $\bi{\{n\}_p}{\{r\}_p}=\bi{0}{t}=0$; if $p\mid r$ then
$$\align\ord_p(T_{l,\al}(n_0,r_0))\gs&\tau_p(\{r_0\}_{p^{\al-1}},\{n_0-r_0\}_{p^{\al-1}})
\\=&\tau_p(p\{r_0\}_{p^{\al-1}},p\{n_0-r_0\}_{p^{\al-1}})
=\tau_p(\{r\}_{p^{\al}},\{n-r\}_{p^{\al}}).
\endalign$$
So we have
$$T_{l,\al+1}(n,r)\eq0\eq(-1)^{\{r\}_p}\bi{\{n\}_p}{\{r\}_p}T_{l,\al}(n_0,r_0)\ (\mo\ p).$$

{\tt Case 3}. $n\eq r\eq0\ (\mo\ p^{\al})$. In this case, $s=t=0$.
When $l=0$, by Theorem 3.1 we have
$$\align &T_{0,\al+1}(n,r)=\f{p^{\lfloor\f{n/p^{\al}-1}{p-1}\rfloor}}{(n/p^{\al})!}
S_{\al+1}^{(p)}\l(\f n{p^{\al}},r\r)
\\\eq&\f{p^{\lfloor\f{n_0/p^{\al-1}-1}{p-1}\rfloor}}{(n_0/p^{\al-1})!}
S_{\al}^{(p)}\l(\f {n_0}{p^{\al-1}},\f rp\r)=T_{0,\al}(n_0,r_0)
\ \l(\mo\ p^{(2-\da_{p,2})(\al-1)}\r).
\endalign$$
As $\al\gs 2$ this implies that $T_{0,\al+1}(n,r)\eq T_{0,\al}(n_0,r_0)\ (\mo\ p)$.
In view of Lemma 2.2 and the induction hypothesis, if $l>0$ then
$$\align T_{l,\al+1}(n,r)=&-\f r{p^{\al}}T_{l-1,\al+1}(n,r+p^{\al+1})
-T_{l-1,\al+1}(n-1,r+p^{\al+1}-1)
\\=&-\f{r_0}{p^{\al-1}}T_{l-1,\al+1}(pn_0,p(r_0+p^{\al}))
\\&-T_{l-1,\al+1}(p(n_0-1)+p-1,p(r_0+p^{\al}-1)+p-1)
\\\eq&-\f{r_0}{p^{\al-1}}T_{l-1,\al}(n_0,r_0+p^{\al})
\\&-(-1)^{p-1}\bi{p-1}{p-1}
T_{l-1,\al}(n_0-1,r_0+p^{\al}-1)
\\\eq&T_{l,\al}(n_0,r_0)\ \ (\mo\ p).
\endalign$$

Combining the above we have completed the induction proof. \qed

\medskip

To establish Theorem 1.8 we need some auxiliary results.

\proclaim{Lemma 4.1} Let $p$ be a prime, and let $n=p^{\beta}(pq+r)$
with $\beta,q,r\in\N$ and $\{-q\}_{p-1}<r<p$. Then
$$\ord_p(n!)=\l\lfloor\f{n-1}{p-1}\r\rfloor\iff q=0.$$
\endproclaim
\Proof. For $s=\{-q\}_{p-1}$ we clearly have
$pq+s\eq0\ (\mo\ p-1)$.
Observe that
$$\align\l\lfloor\f{n-1}{p-1}\r\rfloor=&\f{p^{\beta}-1}{p-1}(pq+r)
+\l\lfloor\f{pq+r-1}{p-1}\r\rfloor
\\=&\f{p^{\beta}-1}{p-1}(pq+r)+\f{pq+s}{p-1}.
\endalign$$
Also,
$$\align \ord_p(n!)=&\sum_{0<i\ls\beta}\f{p^{\beta}(pq+r)}{p^i}
+\sum_{i=1}^{\infty}\l\lfloor\f{p^{\beta}(pq+r)}{p^{\beta+i}}\r\rfloor
\\=&(pq+r)\sum_{0<i\ls\beta}p^{\beta-i}+\ord_p((pq+r)!)
\\=&\f{p^{\beta}-1}{p-1}(pq+r)+\ord_p((pq)!).
\endalign$$

In the case $q=0$,  we have $s=0$ and hence
$$\l\lfloor\f{n-1}{p-1}\r\rfloor=\f{p^{\beta}-1}{p-1}r=\ord_p(n!)$$
by the above.

Now let $q>0$. Then
$$\ord_p((pq)!)=\sum_{i=1}^{\infty}\l\lfloor\f{pq}{p^i}\r\rfloor
<\sum_{i=1}^{\infty}\f{pq}{p^i}=\f{q}{1-p^{-1}}=\f{pq}{p-1}
\ls\f{pq+s}{p-1},$$
and therefore $\ord_p(n!)<\lfloor (n-1)/(p-1)\rfloor$.
This ends the proof. \qed

\Remark\ 4.1. By Lemma 4.1, if $p$ is a prime
and $n$ is a positive integer with $\ord_p(n)=\lfloor\log_pn\rfloor$,
then $\ord_p(n!)=\lfloor(n-1)/(p-1)\rfloor$.
\medskip

Using Lemma 4.1 and Corollary 1.2, we can deduce the following
lemma.

\proclaim{Lemma 4.2} Let $\al\in\N$, $n\in\Z^+$, $r\in\Z$
and $n\eq r\eq0\ (\mo\ 2^{\al})$. Then
$$T_{0,\al+1}^{(2)}(n,r)\eq1\ (\mo\ 2)\iff
 n\ \t{is a power of}\ 2.$$
\endproclaim
\Proof. Clearly
$$T_{0,\al+1}^{(2)}(n,r)=\f{(-1)^r2^{n/2^{\al}-1}}{(n/2^{\al})!}
2^{-\l\lfloor\f{n-2^{\al}}{\varphi(2^{\al+1})}\r\rfloor}
\sum_{k\eq r\,(\mo\ 2^{\al+1})}\bi nk.$$
By Lemma 4.1 in the case $p=2$,
$$\ord_2\l(\f n{2^{\al}}!\r)=\f n{2^{\al}}-1\iff
\f n{2^{\al}}\ \t{is a power of}\ 2.$$

If $\al=0$, then
$$2^{-\l\lfloor\f{n-2^{\al}}{\varphi(2^{\al+1})}\r\rfloor}
\sum_{k\eq r\,(\mo\ 2^{\al+1})}\bi nk=2^{-(n-1)}\sum_{k\eq r\,(\mo\ 2)}\bi nk=1.$$

Now we let $\al\gs1$. Note that both $n_*=n/2^{\al-1}$ and $r_*=r/2^{\al-1}$ are even.
Also, $\{n\}_{2^{\al-1}}=\{r\}_{2^{\al-1}}=0$.
Thus, Corollary 1.2 implies that
$$\align&2^{-\l\lfloor\f{n-2^{\al}}{\varphi(2^{\al+1})}\r\rfloor}
\sum_{k\eq r\,(\mo\ 2^{\al+1})}\bi nk
\eq1\ (\mo\ 2)
\\&\iff n_*=2\ \ \t{or}\ \ n_*\eq2r_*\eq0\ (\mo\ 4).
\endalign$$

Combining the above we find that $T_{0,\al+1}^{(2)}(n,r)\eq1\ (\mo\ 2)$ if and only if
$n$ is a power of $2$. This concludes the proof. \qed

The following result plays a major role in our proof of Theorem 1.8.
\proclaim{Theorem 4.1} Let  $\al,c,d,e\in\N$ with $d<2^e$.
Then
$$T_{l,\al+1}^{(2)}(2^{\al}(2^e+d),2^{\al}c)\eq\da_{l,d}\ (\mo\ 2)
\ \ \t{for all}\ l=0,\ldots,d.\tag4.1$$
\endproclaim
\Proof. By Lemma 4.2, if $n\in\Z^+$, $r\in\Z$ and $n\eq r\eq0\
(\mo\ 2^{\beta})$ with $\beta=1$, then $T_{0,\beta+1}^{(2)}(n,r)\eq
T_{0,\beta}^{(2)}(n/2,r/2)\ (\mo\ 2)$. Thus, by modifying the proof of
Theorem 1.5 (just the third case) slightly, we get a
modified version of (1.7) with $p=2$ and $\al$ replaced by
$\beta=1$. This, together with Theorem 1.5, shows that if
$l\in\N$ and $\al>0$ then
$$\align T_{l,\al+1}^{(2)}(2^{\al}(2^e+d),2^{\al}c)
\eq &T_{l,\al}^{(2)}(2^{\al-1}(2^e+d),2^{\al-1}c)
\\\eq&\cdots\eq T_{l,1}^{(2)}(2^e+d,c)\ \ (\mo\ 2).
\endalign$$
So, it suffices to show the following claim:

 {\tt Claim}. {\it If $l\in\N$, $n\in\Z^+$ and $d_n:=n-2^{\lfloor\log_2 n\rfloor}\gs l$, then
$$T_l(n,r):=T_{l,1}^{(2)}(n,r)\eq\da_{l,d_n}\ (\mo\ 2)\ \ \t{for all}\ r\in\Z.$$}

We use induction on $l$ to show the claim.

As $n\in\Z^+$ is a
power of two if and only if $d_n=0$, in the case $l=0$ the claim follows from Lemma 4.2.

Now let $l\in\Z^+$ and assume the claim for $l-1$.
Let $n\in\Z^+$ with $d_n\gs l$. Clearly $d_1=0<l$ and hence $n>1$.
By Lemma 2.2,
$$T_l(n,r)=-rT_{l-1}(n,r+2)-T_{l-1}(n-1,r+1)\quad\t{for any}\ r\in\Z.$$
Since $d_n>d_{n-1}=d_n-1\gs l-1$, by the induction
hypothesis we have
$$T_{l-1}(n,r+2)\eq\da_{l-1,d_n}\ (\mo\ 2)\ \ \t{and}\ \
T_{l-1}(n-1,r+1)\eq\da_{l-1,d_n-1}\ (\mo\ 2).$$
Therefore
$$T_l(n,r)\eq-r\da_{l-1,d_n}-\da_{l-1,d_n-1}\eq\da_{l,d_n}\ (\mo\ 2).$$
This concludes the induction step, and we are done. \qed

\ms

\noindent{\it Proof of Theorem 1.8.} Write $n=2^{\al}(2^e+d)+c$
with $c,d,e\in\N$, $c<2^{\al}$ and $d<2^e$. Clearly $n_0:=2^e+d=\lfloor n/2^{\al}\rfloor$.
Since $2^e\ls n/2^{\al}<n_0+1\ls 2^{e+1}$, we also have $e=\lfloor\log_2(n/2^{\al})\rfloor$.
By Vandermonde's identity,
$$\align&\sum_{k\eq0\,(\mo\ 2^{\al})}\bi nk(-1)^k\l(\f k{2^{\al}}\r)^l
\\=&\sum_{k\eq0\,(\mo\ 2^{\al})}\sum_{j=0}^c\bi cj
\bi{2^{\al}n_0}{k-j}(-1)^k\l(\f k{2^{\al}}\r)^l
\\=&\sum_{j=0}^c\bi cj(-1)^j\sum_{k\eq-j\,(\mo\ 2^{\al})}
\bi{2^{\al}n_0}k(-1)^k\l(\f{k+j}{2^{\al}}\r)^l.
\endalign$$
If $0<j\ls c<2^{\al}$, then
$$\align&\ord_2\(\sum_{k\eq-j\,(\mo\ 2^{\al})}\bi{2^{\al}n_0}k(-1)^k\l(\f{k+j}{2^{\al}}\r)^l\)
\\&\gs\ord_2(n_0!)+\tau_2(2^{\al}-j,j)>\ord_2(n_0!)
\endalign$$
by the inequality $(1.2)$.
So, we need to show that
$$\sum_{k\eq0\,(\mo\ 2^{\al})}\bi{2^{\al}n_0}k(-1)^k\l(\f k{2^{\al}}\r)^l
=\sum_{k\in\N}\bi{2^{\al}n_0}{2^{\al}k}(-1)^{2^{\al}k}k^l$$
has the same $2$-adic order as $n_0!$. If $k$ is even, then
$\ord_2(k^l)\gs l\gs n_0>\ord_2(n_0!)$.
Thus, it remains to prove that $\ord_2(\Sigma)=\ord_2(n_0!)$,
where
$$\Sigma=\sum_{2\nmid k}\bi{2^{\al}n_0}{2^{\al}k}k^l
=\sum_{j\in\N}\bi{2^{\al}n_0}{2^{\al}(2j+1)}(2j+1)^l.$$

Observe that
$$(2j+1)^l=\sum_{s=0}^l\bi ls2^sj^s
=\sum_{s=0}^l\bi ls2^s\sum_{t=0}^sS(s,t)\ t!\bi jt$$
and hence
$$\align\f{\Sigma}{n_0!}=&\sum_{0\ls t\ls s\ls l}\bi ls 2^sS(s,t)\f{t!}{n_0!}
\sum_{j\in\N}\bi{2^{\al}n_0}{2^{\al+1}j+2^{\al}}\bi jt
\\=&\sum_{0\ls t\ls s\ls l}\bi ls 2^{s-t}S(s,t)(-1)^{2^{\al}}
T_{t,\al+1}^{(2)}(2^{\al}n_0,2^{\al})
\\\eq&\sum_{t=0}^l\bi lt T_{t,\al+1}^{(2)}(2^{\al}n_0,2^{\al})\ (\mo\ 2).
\endalign$$

Now we analyze the parity of $\bi lt T_{t,\al+1}^{(2)}(2^{\al}n_0,2^{\al})$
for each $0\ls t\ls l$.
As $l\gs n_0$ and $l\eq n_0\ (\mo\ 2^e)$, we can write
 $l=n_0+2^eq=2^e(q+1)+d$ with $q\in\N$. Note that
$$\bi ld=\prod_{0<r\ls d}\l(1+\f{2^e}r(q+1)\r)\eq1\ (\mo\ 2).$$
Also, if $0\ls t\ls d$ then $T_{t,\al+1}^{(2)}(2^{\al}n_0,2^{\al})\eq\da_{t,d}\ (\mo\ 2)$
by Theorem 4.1.
When $d<t\ls l$ and $T_{t,\al+1}^{(2)}(2^{\al}n_0,2^{\al})\not=0$, we have
$2^{\al}n_0\gs2^{\al+1}t+2^{\al}$; hence $d<t<n_0/2<2^e$ and thus
$$\align \ord_2\bi lt=&\sum_{i=1}^{\infty}\(\l\lfloor\f l{2^i}\r\rfloor
-\l\lfloor\f t{2^i}\r\rfloor-\l\lfloor\f{l-t}{2^i}\r\rfloor\)
\\\gs&\l\lfloor\f l{2^e}\r\rfloor
-\l\lfloor\f t{2^e}\r\rfloor-\l\lfloor\f{l-t}{2^e}\r\rfloor=(q+1)-0-q>0.
\endalign$$

Combining the above we find that
$$\f{\Sigma}{n_0!}\eq\sum_{t=0}^l\bi ltT_{t,\al+1}^{(2)}(2^{\al}n_0,2^{\al})\eq1\ (\mo\ 2).$$
So $\ord_2(\Sigma)=\ord_2(n_0!)$ as required. \qed

\medskip

\Ack. The joint paper was written during the first author's visit
to the University of California at Irvine as a visiting professor.
He would like to thank Prof. Daqing Wan for the invitation.


\medskip

\widestnumber\key{GKP}

\Refs

\ref\key B\by D. F. Bailey\paper Two $p^3$ variations of Lucas' theorem
\jour J. Number Theory \vol 35\yr 1990\pages 208--215\endref

\ref\key C\by P. Colmez\paper Une correspondance de Langlands locale
$p$-adique pour les representations semi-stables de dimension 2
\finalinfo preprint, 2004\endref

\ref\key DS\by D. M. Davis and Z. W. Sun\paper A number-theoretic
approach to homotopy exponents of {\rm SU}$(n)$ \jour J. Pure Appl. Algebra
\vol 209\yr 2007\pages 57--69\endref

\ref\key D\by  L. E. Dickson\book
History of the Theory of Numbers, {\rm Vol. I}
\publ AMS Chelsea Publ., 1999\endref

\ref\key GKP\by R. L. Graham, D. E. Knuth and O. Patashnik
 \book Concrete Mathematics\publ 2nd ed., Addison-Wesley, New York\yr 1994\endref

\ref\key G\by A. Granville\paper Arithmetic properties of binomial
coefficients.\,I. Binomial coefficients modulo prime powers, {\rm
in: Organic mathematics (Burnaby, BC, 1995), 253--276, CMS Conf.
Proc., 20, Amer. Math. Soc., Providence, RI, 1997}\endref

\ref\key IR\by K. Ireland and M. Rosen
\book A Classical Introduction to Modern Number Theory
{\rm (Graduate texts in math.; 84), 2nd ed.}
\publ Springer, New York, 1990\endref

\ref\key LW\by J.H. van Lint and R. M. Wilson\book A Course in Combinatorics
\publ 2nd ed., Cambridge Univ. Press, Cambridge, 2001\endref

\ref\key S02\by Z. W. Sun\paper On the sum $\sum_{k\eq r\, (\mo\ m)}\bi nk$
and related congruences\jour Israel J. Math.
\vol 128\yr 2002\pages 135--156\endref

\ref\key S03\by Z. W. Sun\paper General congruences for Bernoulli
polynomials\jour Discrete Math.\vol 262\yr 2003\pages 253--276\endref

\ref\key S06\by Z. W. Sun\paper Polynomial extension of Fleck's
congruence\jour Acta Arith.\vol 122\yr 2006\pages 91--100\endref

\ref\key W\by D. Wan\paper Combinatorial congruences and $\psi$-operators
\jour Finite Fields Appl. \vol 12\yr 2006\pages 693--703\endref

\ref\key We\by C. S. Weisman\paper Some congruences for binomial coefficients
\jour Michigan Math. J.\vol 24\yr 1977\pages 141--151\endref

\endRefs
\enddocument